\documentclass[11pt]{article}

\usepackage[margin=1in]{geometry} %1 inch margins
\usepackage{color}
\usepackage{microtype}

\usepackage[plainpages=false, colorlinks=true, allcolors=red]{hyperref}
\usepackage{url}
\usepackage[T1]{fontenc}
\usepackage[utf8]{inputenc}

\usepackage[maxbibnames=99, backend=biber]{biblatex}
\usepackage{breakcites}

\usepackage{amsmath, amsthm, amssymb}
\usepackage{commath} %for \norm and \abs
\usepackage{bm}
\usepackage{dsfont}
\usepackage{mathtools}
\usepackage{mathrsfs}
\usepackage{physics}
\usepackage{stmaryrd}
\allowdisplaybreaks

\usepackage{enumerate} %styles of enumeration
\usepackage{cprotect} %verbatim in caption
\usepackage{fancyvrb}
\usepackage{float}
\usepackage{authblk}

\usepackage[ruled]{algorithm2e}
\newtheoremstyle{algostyle}% name
  {}%      Space above, empty = `usual value'
  {}%      Space below
  {}%         Body font
  {}%         Indent amount (empty = no indent, \parindent = para indent)
  {\bfseries}% Thm head font
  {}%        Punctuation after thm head
  {\newline}% Space after thm head: \newline = linebreak
  {\thmname{#1}\thmnumber{ #2.}}% Thm head spec
\theoremstyle{algostyle}

\usepackage{hyperref}
\usepackage{cleveref}
\makeatletter \@mparswitchfalse \makeatother \normalmarginpar

\usepackage{subcaption}
\usepackage{caption}
\usepackage{empheq}
\usepackage[scr=rsfso]{mathalfa}
\usepackage{euscript}

\usepackage{tikz}
\usepackage{color}
\usetikzlibrary{patterns}

\newcommand{\mb}[1]{\mathbb{#1}}
\newcommand{\mc}[1]{\mathcal{#1}}
\newcommand{\mcb}[1]{\boldsymbol{\mathcal{#1}}}
\newcommand{\ms}[1]{\mathscr{#1}}

\newcommand{\mbu}[1]{\ensuremath{\mathbf{#1}}}

\newcommand{\defeq}{\vcentcolon=}

\newcommand{\mat}[1]{\mathbf{#1}} 
\renewcommand{\vec}[1]{{\mathchoice
            {\mbox{\boldmath$\displaystyle{#1}$}}
            {\mbox{\boldmath$\textstyle{#1}$}}
            {\mbox{\boldmath$\scriptstyle{#1}$}}
            {\mbox{\boldmath$\scriptscriptstyle{#1}$}}}}

\newcommand{\inp}[2]{\langle #1, #2 \rangle}
\newcommand{\E}{\operatorname{\mathbb{E}}}
\newcommand{\Var}{\operatorname{Var}}
\renewcommand{\trace}[1]{\hbox{Tr}#1}
\newcommand{\e}{\varepsilon}
\newcommand{\ntheta}{{n_\theta}}
\newcommand{\muprior}{\mu_{\rm prior}}

\newcommand{\mprior}{m_{\text{prior}}}
\newcommand{\Cprior}{\mc{C}_{\text{prior}}}

\newcommand{\obs}{\mathbf{u}^{\text{obs}}}
\newcommand{\Cnoise}{\mat{\Gamma}_{\text{noise}}}
\newcommand{\Nobs}{N_{\text{obs}}}
\newcommand{\mupost}{\mu_{\rm post}}

\newcommand{\mpost}{m_{{\text{post}}}}
\newcommand{\Cpost}{\mc{C}_{\text{post}}}
\newcommand{\Hmisfit}{\mc{H}_{\rm misfit}}
\newcommand{\ppHmisfit}{\tilde{\mc{H}}_{\rm misfit}}
\newcommand{\KLD}{\mathrm{KLD}}
\newcommand{\infogain}{\Phi_{\rm IG}}
\newcommand{\Einfogain}{\overline{\Phi}_{\rm IG}}
\newcommand{\Stot}{S^{\mathrm{tot}}}
\newcommand{\vtheta}{\vec{\theta}}
\newcommand{\veta}{\vec{\eta}}
\newcommand{\Flin}{\mcb{G}}
\newcommand{\Faff}{\boldsymbol{f}}

\newcommand{\Ns}{N_{s}}

\title{%
  Sensitivity Analysis of the Information Gain in Infinite-Dimensional Bayesian
  Linear Inverse Problems
}
\author{Abhijit Chowdhary}
\author{Alen Alexanderian}
\affil{North Carolina State University, Raleigh, North Carolina, 27607}
\author{Shanyin Tong}
\affil{Columbia University, New York City, New York, 10027}
\author{Georg Stadler}
\affil{New York University, New York City, New York, 10012}

\begin{document}

\maketitle
\abstract{
  \noindent
  We study the sensitivity of infinite-dimensional Bayesian linear inverse problems 
  governed by partial differential equations (PDEs) with respect to modeling
  uncertainties. In particular, we consider derivative-based sensitivity analysis of the
  information gain, as measured by the Kullback--Leibler divergence from the posterior
  to the prior distribution.  To facilitate this, we develop a fast and accurate method
  for computing derivatives of the information gain with respect to auxiliary model
  parameters.  Our approach combines low-rank approximations, adjoint-based eigenvalue
  sensitivity analysis, and post-optimal sensitivity analysis.  The proposed approach
  also paves way for global sensitivity analysis by computing derivative-based global
  sensitivity measures.  We illustrate different aspects of the proposed approach using
  an inverse problem governed by a scalar linear elliptic PDE, and an inverse problem
  governed by the three-dimensional equations of linear elasticity, which is motivated
  by the inversion of the fault-slip field after an earthquake.
} 

\section{Introduction}
We consider Bayesian inverse problems governed by partial differential
equations (PDEs) to estimate parameters and their uncertainty. In particular, we
focus on linear inverse problems, where we consider estimation of an \emph{inversion
  parameter} $m$ (usually a function) from the model
\[
  \obs = \mcb{F}m + \veta.
\]
Here, $\obs$ is a vector of measurement data, $\mcb{F}$ is a linear
parameter-to-observable map, and $\veta$ models measurement error.  We consider the case
where computing $\mcb{F}m$ requires solving a system of PDEs.  In such problems, the
governing PDEs typically have parameters that are not or cannot be estimated, but are
required for a complete model specification.  We call such uncertain parameters
\textit{auxiliary parameters}.  Solving such parameterized inverse problems is common in
applications.  For example, an inverse problem may be formulated to estimate an unknown
source term in a PDE model, with parameters in the boundary conditions or coefficient
functions that are also uncertain. In practice, not all uncertain model parameters can
be estimated simultaneously. Reasons for this include limitations in computational
budget, lack of data that is informative to all the parameters, or aleatoric
(irreducible) uncertainty in auxiliary parameters.  Let $\vtheta$ denote the vector of
auxiliary model parameters.  To emphasize the dependence of the parameter-to-observable
map on $\vtheta$, we denote $\mcb{F} = \mcb{F}(\vtheta)$.  In the present work, we seek
to address the following question: how sensitive is the solution of the Bayesian inverse
problem to the different components of $\vtheta$?

In a Bayesian inverse problem, we use prior knowledge, encoded in a prior distribution
law $\muprior$, a model, and measurement data $\obs$ to obtain the posterior
distribution $\mupost$ of $m$.  The posterior is a distribution law for $m$ that is
consistent with the prior and the observed data. Note, however, that the auxiliary
parameters might have a significant influence on the solution of the Bayesian inverse
problem. When these auxiliary parameters are uncertain, this can pose a problem.
Therefore, it is important to have a systematic way of assessing the sensitivity of
$\mupost$ to $\vtheta$.  This can provide important insight and can help determine which
auxiliary parameters require extra care in their specification.  Also, if one has access
to data that informs the important auxiliary parameters, the inverse problem may be
reformulated to include such parameters in the inversion process.

Performing sensitivity analysis on the posterior measure itself is challenging.  In the
present work, we focus on sensitivity of a specific computable quantity of interest
(QoI) defined in terms of $\mupost$.  Namely, we consider the sensitivity analysis of
the Kullback--Leibler divergence (KLD) \cite{Kullback_Leibler_1951,Gibbs_Su_2002} of the
posterior from the prior:
\begin{equation} \label{eq:information-gain}
  \KLD(\mupost || \muprior)
  \defeq \int \log \left[ \dv{\mupost}{\muprior} \right] \dd{\mupost}.
\end{equation}
This quantity, which is also known as the relative entropy of the posterior with respect
to the prior, provides a measure of \emph{information gain} in the process of Bayesian
inversion; see \Cref{sec:math-prelim-information-gain}.  In this article, we focus on
the construction of scalable methods for the sensitivity analysis of this QoI, for
linear inverse problems governed by PDEs with infinite-dimensional inversion parameters.

\subsection{Related work}
The present work builds on the efforts in hyper-differential sensitivity
analysis (HDSA) or post-optimality sensitivity analysis
\cite{Sunseri_Hart_Waanders_Alexanderian_2020,
Sunseri_Alexanderian_Hart_Waanders_2022,Waanders_Hart_Herzog_2020,
Brandes_Griesse_2006,Griesse_2007}.  Originally, HDSA was intended to enable differentiating
the solution of a deterministic optimization problem with respect to auxiliary
parameters in the objective function.  Further work, specifically
\cite{Sunseri_Alexanderian_Hart_Waanders_2022}, considered such an analysis for
the maximum a posteriori probability (MAP) point in Bayesian inverse problems
and the Bayes risk as the HDSA QoI.  We seek to take a step further towards QoIs
that incorporate more information than the MAP point.

There are various options for QoIs that measure different aspects of posterior
uncertainty. For example, we can consider Bayesian optimal experimental design (OED)
objectives such as the Bayesian A- and D-optimality
criteria~\cite{Chaloner_Verdinelli_95,Alexanderian_Saibaba_2018,Alexanderian_21}.
Notably, the D-optimality design criterion, otherwise known as the expected information
gain, is the expectation of the information gain over data. While the D-optimality
criterion is suitable for OED, as discussed in this article, the information gain itself
provides important additional insight. Both the information gain and its expectation
over data admit closed-form expressions for Gaussian linear inverse problems; for
infinite dimensional problems, they are given in~\cite{Alexanderian_Gloor_Ghattas_2016}.
These expressions are important for this present work.

\subsection{Our approach and contributions}

In the present work, we focus on HDSA of the information gain. This is illustrated
briefly in~\Cref{sec:simple} with a simple example.  Computing
\eqref{eq:information-gain} in large-scale inverse problems is challenging. Even with
the assumption of a Gaussian linear inverse problem for which
\eqref{eq:information-gain} admits a closed-form expression, it still involves
infinite-dimensional operators as well as the MAP point.  Derivative-based sensitivity
analysis of this QoI adds an additional layer of complexity.

The main contributions of this work are in developing a scalable computational method
for HDSA of the information gain and the expected information gain in
infinite-dimensional Bayesian linear inverse problems governed by PDEs.  Our approach,
presented in \Cref{sec:methods}, combines fast methods for post-optimal sensitivity
analysis of the MAP point, low-rank approximations, and adjoint-based eigenvalue
sensitivity analysis for operators involved in $\KLD(\mupost || \muprior)$. Our
developments include a detailed computational algorithm and a discussion of its
computational cost.  The proposed approach also makes derivative-based \emph{global}
sensitivity analysis (see~\Cref{sec:math-prelim-sobol}) feasible.

We present comprehensive numerical experiments that demonstrate various aspects of the
proposed approach (see \Cref{sec:numerics}).  We consider two PDE-constrained inverse
problems.  The first one, discussed in \Cref{sec:elliptic-toy}, illustrates the
difference between HDSA of the information gain versus HDSA of the \emph{expected}
information, in an inverse problem governed by a (scalar) linear elliptic PDE.  The
second example is an inverse problem governed by equations of linear elasticity and it
is motivated by applications in fault-slip reconstruction in seismology. We fully
elaborate our proposed approach for that example; see \Cref{sec:cr-model}.

\subsection{An illustrative example}\label{sec:simple}
We use a simple analytic example from~\cite{Youssef_Feng_2019} to illustrate the
information gain as a measure of posterior uncertainty.  The inverse problem is to infer
$\mbu{m} \in \mathbb R^2$ in the model 
\begin{equation}
  \obs = \mat{F}(\vtheta) \mbu{m} + \veta,
  \quad \text{ where } \quad
  \mat{F}(\vtheta) =
  \begin{bmatrix}
    \theta_2 & \theta_1 \\ \theta_1 & 1 - \theta_2
  \end{bmatrix}.
\end{equation}
Here $\obs$ is a vector of observations and $\veta \sim \mc{N}(\vec{0}, \sigma^2
\mat{I})$ represents additive Gaussian errors.  In this problem, $\vtheta = [\theta_1\;
\theta_2]^\top \in [0,1]^2$ is the vector of auxiliary parameters. Additionally, we
assume prior knowledge of $\mbu{m}$ is encoded in a Gaussian prior $\muprior =
\mc{N}(\vec{0}, \mat{I})$. Due to linearity of the parameter-to-observable map and use
of Gaussian prior and noise models, the posterior is also a Gaussian given by
\begin{equation}
  \mupost
  \sim
  \mc{N}\left(
  \mbu{m}_{\rm post}, \mat{C}_{\rm post}
  \right),
  \quad \text{ where } \quad
  \begin{cases}
    \mbu{m}_{\rm post} = \sigma^{-2} \mat{C}_{\rm post} \mat{F} \obs, \\
    \mat{C}_{\rm post} = (\sigma^{-2}\mat{F}^\top \mat{F} + \mat{I})^{-1}.
  \end{cases}
\end{equation}
In the present (finite-dimensional) Gaussian linear setting, the information gain takes
the following expression:
\begin{equation}\label{eq:Phi_KL}
  \KLD(\mupost || \muprior)
  = \frac{1}{2}\left[
    \log\det(\mat{C}_{\rm post}^{-1})
    + \trace(\mat{C}_{\rm post})
    - 2
    + \inp{\mbu{m}_{\rm post}}{\mbu{m}_{\rm post}}
    \right].
\end{equation}
\begin{figure}[!htb]
  \centering
  \begin{subfigure}[b]{0.495\textwidth}
    \caption{
      $|\partial_{\theta_1} \infogain|$ \vspace{.2em}
    }
    \includegraphics[width=\textwidth]{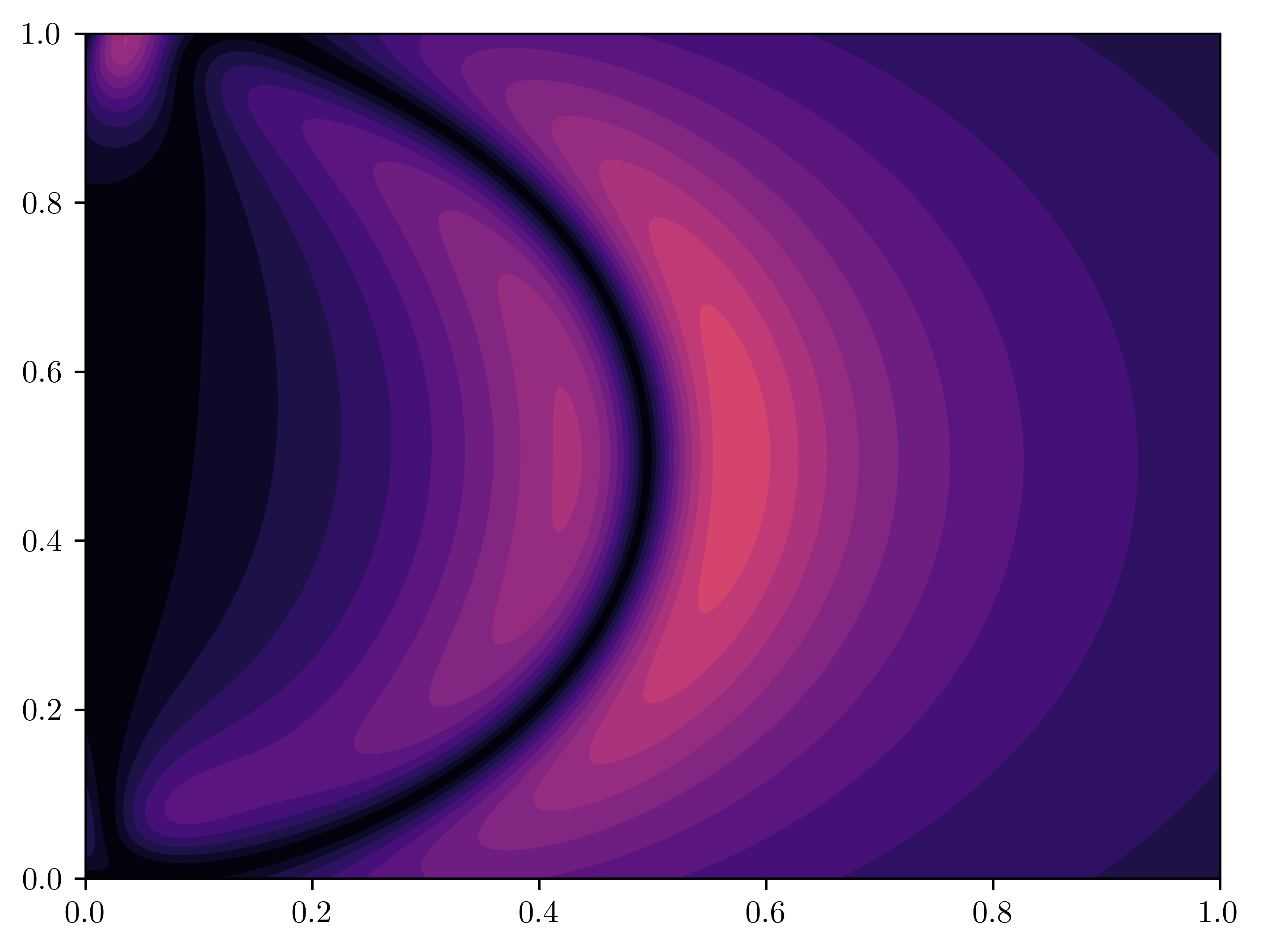}
  \end{subfigure}
  \begin{subfigure}[b]{0.495\textwidth}
    \caption{
      $|\partial_{\theta_2} \infogain|$ \vspace{.2em}
    }
    \includegraphics[width=\textwidth]{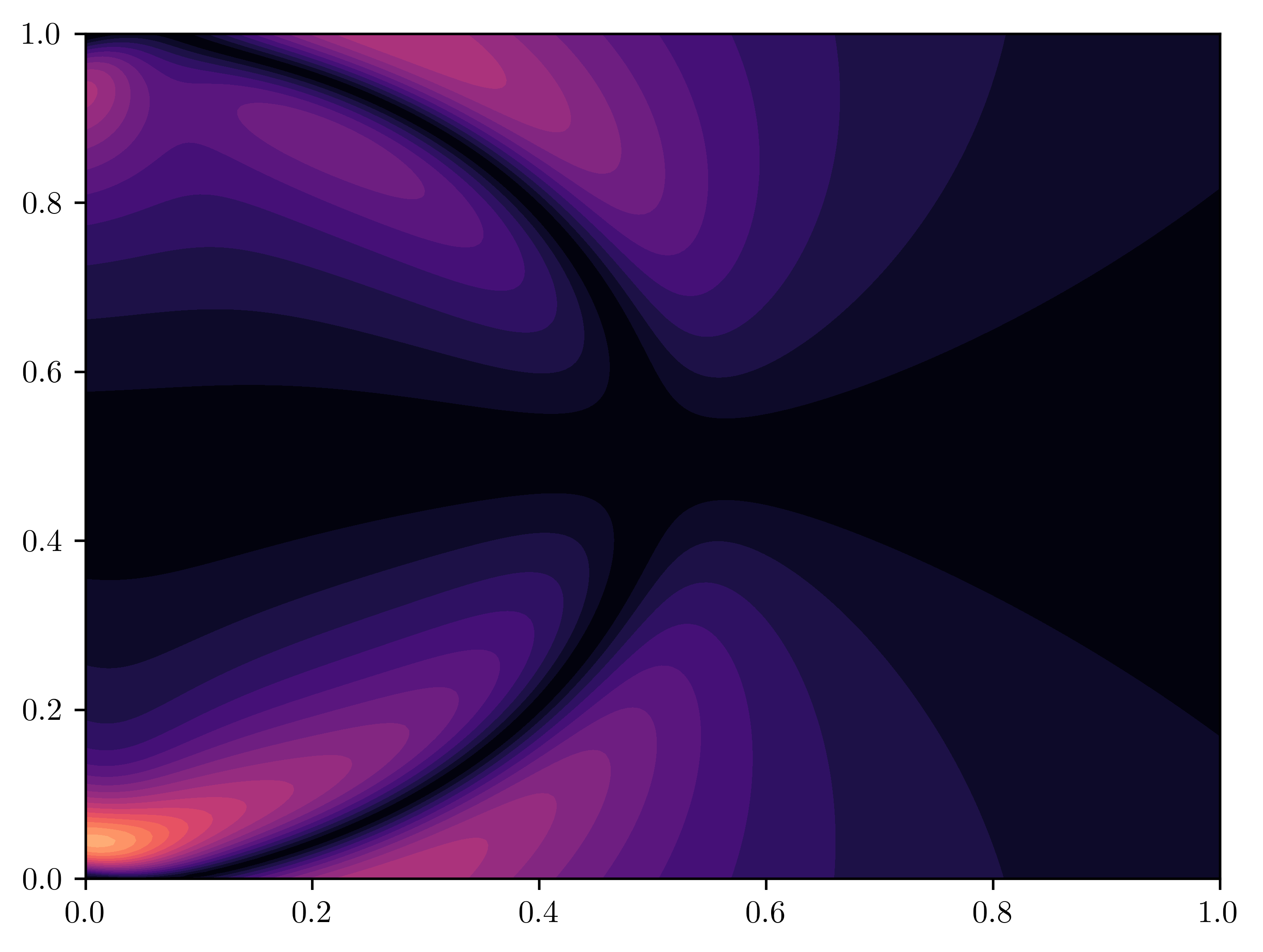}
  \end{subfigure}
  \includegraphics[width=0.500\textwidth]{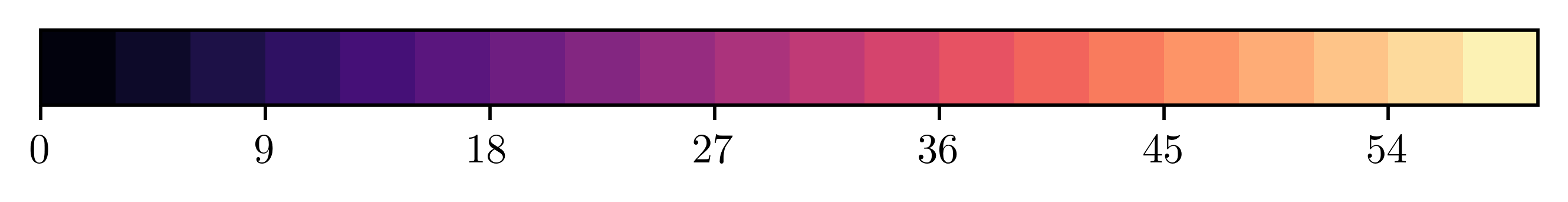}
  \caption{
    Information gain sensitivities for two-dimensional model problem. Shown are the
    absolute values of the derivatives with respect to the auxiliary parameters.
  } \label{fig:simple}
\end{figure}

The main target of this paper is to compute and understand derivatives of $\KLD$ with
respect to auxiliary parameters such as $\vtheta$. In this example, these derivatives
are straightforward to compute and their absolute values are shown in \Cref{fig:simple}
for $\obs=[0.15\; 0.05]^\top$ and $\sigma=0.1$. Note that even for this simple example,
the sensitivities exhibit complex behavior and change substantially for small changes in
the auxiliary parameters. This is a common feature of sensitivity analysis in inverse
problems.

The expression \eqref{eq:Phi_KL} for the information gain already hints towards
challenges for infinite or high-dimensional parameters. For instance, the trace and
determinant of the posterior covariance operator may be challenging to approximate, and
$\mat{C}_{\rm post}$ and $\mbu{m}_{\rm post}$ must be differentiated with respect to the
auxiliary parameters $\vtheta$.  We will show how such a sensitivity analysis can be
done for PDE-constrained linear inverse problems with infinite-dimensional parameters.

\section{Preliminaries}
\label{sec:math-prelim}
In this section, we provide the background material for Bayesian linear inverse
problems, adjoint-based gradient and Hessian computation, the information gain, and
global sensitivity analysis.

\subsection{Linear Bayesian Inverse Problems}
Assume the forward model is governed by a PDE with its weak formulation represented
abstractly as follows: the solution $u \in \mc{U}$ satisfies
\begin{equation} \label{eq:abstract-pde}
  \ms{A}(p, u;\vtheta)
  + \ms{C}(p, m; \vtheta)
  + \ms{D}(p; \vtheta)
  = 0
  , \qquad \text{for all } p \in \mc{V}.
\end{equation}
Here, $\mc{U}$ and $\mc{V}$ are appropriately chosen infinite-dimensional Hilbert
spaces.  The inversion parameter $m$ belongs to an infinite-dimensional Hilbert space
$\mc{M}$.  The set of the auxiliary parameters is denoted as $\Theta$.  For simplicity,
we assume this set is of the form $\Theta \defeq \Theta_1 \times \dots \times
\Theta_{\ntheta}$ with $\Theta_i \subset \mb{R}$ for each $i \in \{1, \ldots,
\ntheta\}$.  The form \eqref{eq:abstract-pde} is chosen to mimic a common type of weak
form arising from linear PDEs that are linear in the inversion parameter and potentially
contain non-homogeneous volume or boundary source terms. Therefore, we assume that all
forms are linear in all inputs except for $\vtheta$.

Typically, $\ms{A}$ represents a differential operator on the state variable, $\ms{C}$
is the term involving the inversion parameter, and $\ms{D}$ aggregates any miscellaneous
boundary conditions or source terms that do not depend on the inversion parameters. As
an example, consider the following elliptic PDE
\begin{equation}
  \begin{alignedat}{2}
    -\Delta u + \theta u  &= m &\quad \text{in } &\Omega , \\
    \nabla u \cdot n &= g &\quad \text{on } &\partial\Omega,
  \end{alignedat}
\end{equation}
with the source term $m$ as the inversion parameter, and given boundary data $g$ and
unit boundary normal $n$.  This problem has the weak form \eqref{eq:abstract-pde} with
\begin{equation}
  \ms{A}(p,u; \theta) = \int_\Omega \nabla u \cdot \nabla p + \theta u p \dd{x},
  \quad
  \ms{C}(p, m) = \int_\Omega mp \dd{x},
  \quad
  \ms{D}(p) = \int_{\partial\Omega} gp \dd{s},
\end{equation}
along with $\mc{U} = \mc{V} = H^1(\Omega)$ and $\mc{M} = L^2(\Omega)$. Note, the
auxiliary parameter here only appears in $\ms{A}$, but in general it can appear in any
of the forms.

Next, we consider a Bayesian inverse problem governed by \eqref{eq:abstract-pde}.  We
assume that we have a vector of measurement data $\obs$ that relates to the inversion
parameter $m$ via
\begin{equation}
  \obs = \mcb{F}(\vtheta)(m) + \veta.
\end{equation}
Here, $\mcb{F}$ is the parameter-to-observable map whose evaluation typically breaks
down into two steps: solving \eqref{eq:abstract-pde} for $u$ and evaluating $u$ at the
measurement points.  The vector $\veta \in \mb{R}^{\Nobs}$ represents measurement error.
We assume $\veta \sim \mc{N}(\vec{0}, \Cnoise)$.  Additionally, we assume $\muprior =
\mc{N}(\mprior, \Cprior)$.  The specific choice of $\Cprior$ is important.  This
operator must be strictly positive, self-adjoint, and trace-class.  A convenient
approach for constructing such a $\Cprior$ is to define it as the inverse of
Laplacian-like differential operator; see \cite{Stuart_2010} for details.  The prior
measure induces the Cameron--Martin space $\mc{E} = {\rm range}(\Cprior^{1/2})$ which is
equipped with the inner product
\begin{equation}
  \inp{x}{y}_{\Cprior^{-1}} \defeq \inp{\Cprior^{-1/2}x}{\Cprior^{-1/2}y}
  , \qquad \text{for all } x, y \in \mc{E}.
\end{equation}
Finally, we further assume that $\mprior \in \mc{E}$.

With the data model and prior measure defined, Bayes' formula in infinite dimensions
\cite{Bui-Thanh_Ghattas_Martin_Stadler_2013} states that
\begin{equation}
  \dv{\mu_{\rm post}}{\muprior} \propto \pi_{\rm like}(\obs | m; \vtheta),
\end{equation}
where $\pi_{\rm like}$ is the likelihood probability density function 
\begin{equation}
  \pi_{\rm like}(\obs | m; \vtheta) \propto \exp \left[
  -\frac{1}{2} \| \mcb{F}(\vtheta)(m) - \obs \|_{\Cnoise^{-1}}^2
  \right].
\end{equation}
For a linear parameter-to-observable map, using Gaussian prior and noise models implies
that the posterior measure is also Gaussian and has a simple closed
form~\cite{Stuart_2010}.  Note that in our abstract weak form~\eqref{eq:abstract-pde},
the presence of nonzero $\ms{D}$ makes $\mcb{F}$ affine as opposed to linear. This
results in a Gaussian posterior with very similar closed form to the linear case, albeit
with a slight change in the expression for the MAP point.  This can, however, have a
significant impact on the solution of the inverse problem.  Henceforth, assuming an
affine $\mcb{F}$, we write its action as
\begin{equation}\label{equ:parameter_to_observable}
  \mcb{F}(m) = \Flin m + \Faff.
\end{equation}
Here, $\Flin : \mc{M} \to \mb{R}^{\Nobs}$ is a continuous linear transformation and
$\Faff \in \mb{R}^{\Nobs}$. The resulting closed form expressions for the mean and
covariance of the Gaussian posterior are \cite{Stuart_2010}
\begin{equation} \label{eq:closed-form-posterior}
  \mu_{\rm post} = \mc{N}(\mpost, \Cpost),
  \quad \text{where} \quad
  \begin{cases}
    \Cpost = (\Flin^* \Cnoise^{-1} \Flin + \Cprior^{-1})^{-1} , \\
    \mpost = \Cpost \left[
      \Flin^* \Cnoise^{-1} (\obs - \Faff) + \Cprior^{-1} \mprior
      \right].
  \end{cases}
\end{equation}

\subsection{Adjoint-Based Gradient and Hessian Computation}
Generally, the MAP point $\mpost$ can be obtained by minimizing
\begin{equation} \label{eq:objective-functional}
  \mc{J}(m; \vtheta)
  \defeq
  \frac{1}{2} \| \mcb{F}(\vtheta)(m) - \obs \|_{\Cnoise^{-1}}^2
  + \frac{1}{2} \| m - \mprior \|_{\Cprior^{-1}}^2,
\end{equation}
over $\mc{E}$.
In the linear Gaussian settings, this is a ``regularized'' linear least squares problem.
In what follows, we need the adjoint-based
expressions~\cite{Gunzburger_2002,Petra_Stadler_2011,Ghattas_Willcox_2021} for the
gradient and Hessian of $\mc{J}$ with respect to $m$.  These expressions and the
approach for their derivation are central to our proposed methods.  Assume that
\eqref{eq:abstract-pde} has a unique solution $u(m; \vtheta)$ for every $m \in \mc{M}$
and $\vtheta\in \Theta $. Next, we define a linear observation operator $\mcb{Q}:
\mc{V} \to \mb{R}^{\Nobs}$, which extracts solution values at the measurement points.
Thus, we can write $\mcb{F}(\vtheta)(m) = \mcb{Q}u(m; \vtheta)$.

To facilitate derivative computation, we follow a formal Lagrange approach.  To this
end, we consider the Lagrangian
\begin{equation}
  \mc{L}(u,p,m)
  \defeq
  \frac{1}{2} \| \mcb{Q}u - \obs \|_{\Cnoise^{-1}}^2
  +\frac{1}{2} \| m - \mprior \|_{\Cprior^{-1}}^2
  + \ms{A}(p, u; \vtheta)
  + \ms{C}(p, m; \vtheta)
  + \ms{D}(p; \vtheta).
\end{equation}
The gradient of $\mc{J}$ at $m$ in direction $\tilde{m}$, denoted by
$\mc{G}(m)(\tilde{m})$, is given by 
\[
   \mc{G}(m)(\tilde{m}) \defeq \mc{L}_m(u,p,m)(\tilde{m}),             
\]
where $u$ and $p$, respectively, satisfy the state and adjoint equations:
\begin{equation*}
  \begin{aligned}
    \mc{L}_p(u,p,m)(\tilde{p}) &= 0, \quad\forall \tilde{p} \in \mc{V}, \\
    \mc{L}_u(u,p,m)(\tilde{u}) &= 0, \quad\forall \tilde{u} \in \mc{V}. \\
  \end{aligned}
\end{equation*}
For the class of inverse problems under study, these 
expressions take the form  
\begin{equation} \label{eq:abstract-gradient-system}
  \begin{alignedat}{2}
    &\mc{G}(m)(\tilde{m})
    = \inp{m-\mprior}{\tilde{m}}_{\Cprior^{-1}}
    + \ms{C}(p,\tilde{m};\vtheta), & \\
    &\text{where}&
  \\
    &\ms{A}(\tilde{p},u;\vtheta)
    + \ms{C}(\tilde{p},m;\vtheta)
    + \ms{D}(\tilde{p}; \vtheta)
    = 0,
    &\quad \forall \tilde{p} \in \mc{V}, \\
    &\inp{\mcb{Q}u-\obs}{\mcb{Q}\tilde{u}}_{\Cnoise^{-1}}
    + \ms{A}(p, \tilde{u};\vtheta)
    = 0,
    &\quad \forall \tilde{u} \in \mc{V}.
  \end{alignedat}
\end{equation}
Likewise, we can define the Hessian of \eqref{eq:objective-functional} by considering
the meta-Lagrangian~\cite{Ghattas_Willcox_2021}:
\begin{equation} \label{eq:lagrangian-hessian}
  \mc{L}^H(u,p,m,\hat{u},\hat{p},\hat{m})
  \defeq \mc{L}_m(\hat{m}) + \mc{L}_u(\hat{u}) + \mc{L}_p(\hat{p}).
\end{equation}
Note that, for notational convenience, we have suppressed the dependence of $\mc{L}$ on
$(u,p,m)$.  The Hessian action $\mc{H}(\hat{m},\tilde{m})$ is given by
$\mc{H}(\hat{m},\tilde{m}) \defeq
\mc{L}^H_{m}(u,p,m,\hat{u},\hat{p},\hat{m})(\tilde{m})$, where $\hat{u}$ and $\hat{p}$
are obtained, respectively, by setting the variations of $\mc{L}^H$ with respect to $p$
and $u$ to zero. For the case of linear inverse problems under study,
\begin{subequations} \label{eq:abstract-incremental-system}
  \begin{align}
    \label{eq:abstract-hessian-action}
    &	\mc{H}(\hat{m},\tilde{m})
    = \inp{\tilde{m}}{\hat{m}}_{\Cprior^{-1}}
    + \ms{C}(\hat{p},\tilde{m}; \vtheta), &
    \\
    &\text{where}\notag\\
    \label{eq:abstract-incremental-state}
    &	\ms{A}(\tilde{p},\hat{u};\vtheta)
    + \ms{C}(\tilde{p},\hat{m};\vtheta)
    = 0,
    \quad\forall \tilde{p} \in \mc{V}, \\
    \label{eq:abstract-incremental-adjoint}
    &	\inp{\mcb{Q}\tilde{u}}{\mcb{Q}\hat{u}}_{\Cnoise^{-1}}
    + \ms{A}(\hat{p},\tilde{u};\vtheta)
    = 0,
    \quad\forall \tilde{u} \in \mc{V}. 
  \end{align}
\end{subequations}
Note that, since we focus on linear/affine inverse problems, the Hessian has no
dependency on the inversion parameter $m$.  The
equations~\eqref{eq:abstract-incremental-state} and
\eqref{eq:abstract-incremental-adjoint} are the so-called incremental state and
incremental adjoint equations.  Of additional interest is the Hessian of only the data
misfit term in \eqref{eq:objective-functional}, denoted as the data-misfit Hessian
$\Hmisfit$.  Following a similar procedure to the above, we have 
\begin{equation} \label{eq:abstract-Hmisfit-system}
  \begin{alignedat}{2}
    &\Hmisfit(\hat{m},\tilde{m})
    = \ms{C}(\hat{p},\tilde{m};\vtheta), &\\
    &\text{where}& \\
    &\ms{A}(\tilde{p},\hat{u};\vtheta)
    + \ms{C}(\tilde{p},\hat{m};\vtheta)
    = 0, &\quad\forall \tilde{p} \in \mc{V}, \\
    &\inp{\mcb{Q}\tilde{u}}{\mcb{Q}\hat{u}}_{\Cnoise^{-1}}
    + \ms{A}(\hat{p},\tilde{u};\vtheta)
    = 0, &\quad\forall \tilde{u} \in \mc{V}.
  \end{alignedat}
\end{equation}

\subsection{Quantities of Interest for Bayesian Inference: The Information Gain}
\label{sec:math-prelim-information-gain}
To perform sensitivity analysis on the posterior distribution of an inverse problem, we
consider scalar quantities of interest (QoIs) that describe specific aspects of the
posterior measure.  As mentioned in the introduction, one possibility is to consider
optimal experimental design (OED) criteria, which provide different measures of
uncertainty in the parameters.  A classical approach to deriving Bayesian OED criteria
is to adapt ideas from information theory.  In that context, a natural metric for
individual distributions is (Shannon) entropy.  However, the continuous analogue of
entropy, the differential entropy, lacks many fundamental properties of its discrete
counterpart.  Instead, it is common to use the relative entropy, as given by the
\textit{Kullback--Leibler divergence (KLD)} from the posterior to the prior:
\begin{equation*}
  \KLD(\mupost || \muprior)
  \defeq
  \int_\mc{M} \log \left[ \dv{\mupost}{\muprior} \right] \dd{\mupost}.
\end{equation*}
The notion of relative entropy has seen much recent popularity in the machine learning
literature, notably for variational auto-encoders \cite{Kingma_Welling_2013} and
information bottlenecks \cite{Tishby_Zaslavsky_2015}. In Bayesian analysis, the relative
entropy provides a measure of information gain. In that context, the expected relative
entropy, also known as the \textit{expected information gain}, is a classical choice for
a design criterion. This is also known as the Bayesian D-optimality criterion.  This
expected information gain is defined as the expectation of the KLD over data:
\begin{equation}
  \label{eq:expected-information-gain}
  \overline{\KLD}(\mupost || \muprior)
  \defeq \E_{\obs} \left[ \KLD(\mupost || \muprior) \right].
\end{equation}
For linear Gaussian inverse problems, the information gain has a closed-form expression.
Furthermore, it has been shown that a similar analytic expression can be derived in the
infinite-dimensional setting~\cite{Alexanderian_Gloor_Ghattas_2016}:
\begin{equation}
  \label{eq:linear-kld}
  \KLD(\mupost || \muprior)
  \defeq
  \frac{1}{2} \left[
  \log \det ( \ppHmisfit + I )
  - \trace(\Hmisfit\Cpost)
  + \inp{\mpost - \mprior}{\mpost - \mprior}_{\Cprior^{-1}}
  \right],
\end{equation}
where $\ppHmisfit$ is the prior-preconditioned data-misfit Hessian, given by
\begin{equation}
  \label{eq:prior-preconditioned-data-misfit-hessian}
  \ppHmisfit
  \defeq
  \Cprior^{1/2} \Hmisfit \Cprior^{1/2}.
\end{equation}
We seek to assess the (derivative-based) sensitivity of quantities such as the
information gain to the auxiliary parameters in the governing PDEs.  As before, we let
$\vtheta$ denote a vector of auxiliary parameters. In this setting, the posterior
measure depends on $\vtheta$ as well.  In what follows, we denote $\infogain(\vtheta)
\defeq \KLD(\mupost(\vtheta) || \muprior)$.  Likewise, we write the expected information
gain as $\Einfogain(\vtheta) \defeq \overline{\KLD}(\mupost(\vtheta) || \muprior)$. Note
that for linear Gaussian inverse problems,
\[
  \Einfogain(\vtheta) = \frac12 \log \det ( \ppHmisfit + I ).
\]

\subsection{Global Sensitivity via Derivative-Based Upper Bounds for Sobol Indices}
\label{sec:math-prelim-sobol}
A standard derivative-based sensitivity analysis of the information gain provides only
local sensitivity information.  It is also of interest to perform global sensitivity
analysis.  A standard approach in the uncertainty quantification literature is to use a
variance-based global sensitivity analysis approach. This involves computing
variance-based indices known as Sobol indices~\cite{Sobol_1990,Sobol_2001}. Computing
such indices requires an expensive sampling procedure, which would be computationally
challenging for the quantities of interests considered in the present work.  This is due
to the need for a potentially large number of QoI evaluations.  Instead, we rely on
derivative-based upper bounds for total Sobol indices
\cite{SobolKucherenko09,Kucherenko_2016}, which provide a tractable approach for global
sensitivity analysis of the information gain.

We assume that the entries of $\vtheta \in \Theta$ are independent random variables with
cumulative density functions (CDFs) $F_1, \dots, F_{n_\theta}$. Also, we assume
$\infogain \in L^2(\Theta, \sigma(\vtheta), F)$, where $\sigma(\vtheta)$ is the
$\sigma$-algebra generated by $\{\theta_i\}_{i=1}^{n_\theta}$, and $F$ is the joint CDF
of $\vtheta$. The first order Sobol indices are given by
\begin{equation} \label{eq:first-order-sensitivity-index}
  S_i \defeq \frac{
    \Var[\E[\infogain(\vtheta) | \theta_i]]
  }{
    \Var[\infogain(\vtheta)]
  }, \; i = 1, \dots, n_\theta.
\end{equation}
Here, $S_i$ can be interpreted as quantifying the contribution of $\theta_i$ on
$\Var[\infogain(\vtheta)]$. One can also define higher order sensitivity indices that,
in addition to first order effects, account for interactions between random inputs.  To
this end, we consider the total Sobol indices
\begin{equation} \label{eq:total-sensitivity-index}
  \Stot_i \defeq \frac{
    \E[\Var[\infogain(\vtheta) | \vtheta^{\sim i}]]
  }{
    \Var[\infogain(\vtheta)]
  }
  = \frac{
    \Var[\infogain(\vtheta)]
    - \Var[\E[\infogain(\vtheta) | \vtheta^{\sim i}]]
  }{
    \Var[\infogain(\vtheta)]
  },
\end{equation}
where
\begin{equation}
  \vtheta^{\sim i}
  \defeq [ \theta_1 \; \cdots \; \theta_{i-1} \; \theta_{i+1} \; \cdots, \theta_{n_\theta}]^\top
  \in \mb{R}^{n_\theta - 1}.
\end{equation}
Note that $\Stot_i$ quantifies the contribution of $\theta_i$ by itself, and through its
interactions with the other entries of $\vtheta$, to $\Var[\infogain(\vtheta)]$.

It has been shown in \cite{Kucherenko_2016}, utilizing the $L^2$-Poincar\'{e} inequality
\begin{equation}\label{eq:poincare}
  \int \infogain(\vtheta)^2 \dd{F(\vtheta)}
  \leq C(F) \int \|\nabla \infogain(\vtheta)\|^2 \dd{F(\vtheta)},
\end{equation}
that these total sensitivity indices obey the upper bound
\begin{equation}\label{eq:sobol-bound}
  \Stot_i \leq \frac{C(F_i)}{\Var[\infogain(\vtheta)]}
  \int \left(
  \pdv{}{\theta_i}\infogain(\vtheta)
  \right)^2 \dd{F(\vtheta)}.
\end{equation}
Here, $C(F_i)$ is the Poincare constant corresponding to CDF $F_i$.  For example, if
$\theta_i$ is uniformly distributed random variable on $[-1, 1]$, the corresponding
Poincar\'{e} constant is $C(F_i) = 4 / \pi^2$;
see~\cite{SobolKucherenko09,Kucherenko_2016}.

The above discussion shows that an efficient computational approach for computing
derivatives of the information gain at arbitrary choices of $\vtheta$ enables the
approximation of an upper bound for the total order Sobol indices.  This provides global
sensitivity measures  that take the ranges of uncertainty in the auxiliary parameters
into account.

\section{Method} \label{sec:methods}
In this section, we outline a scalable computational approach for sensitivity analysis
of \eqref{eq:linear-kld} with respect to the auxiliary parameters.  Our approach
combines a low-rank spectral decomposition (see~\Cref{sec:low-rank}) and adjoint-based
eigenvalue sensitivity analysis (see~\Cref{sec:eigen-sens}).  We also need to borrow
tools from post-optimal sensitivity analysis to compute the derivative of the term
involving the MAP point, in the definition of~\eqref{eq:linear-kld}. This is outlined
in~\Cref{sec:HDSA}.  In~\Cref{sec:alg}, we summarize the overall computational procedure
for sensitivity analysis of the information gain, and discuss the computational cost of
the proposed approach, in terms of the required number of PDE solves.

\subsection{Low-Rank Approximation of the Information Gain}\label{sec:low-rank}
Consider the data-misfit Hessian described implicitly in
\eqref{eq:abstract-Hmisfit-system}, which has the following closed form for affine
inverse problems under study
\begin{equation}
  \Hmisfit \defeq \Flin^* \Cnoise^{-1} \Flin,
\end{equation}
where $\Flin$ is as in~\eqref{equ:parameter_to_observable}.  This is a positive
self-adjoint trace-class operator, and thus has a spectral decomposition with
non-negative eigenvalues and orthonormal eigenvectors.  Recall also that the
observations are finite-dimensional and belong to $\mb{R}^{\Nobs}$. Thus, the range of
$\Hmisfit$ has finite dimension $\Nobs$.  Moreover, in inverse problems governed by PDEs
with smoothing forward operators the eigenvalues of this Hessian operator decay rapidly
and the numerical rank of $\Hmisfit$ is typically smaller than $\Nobs$.  In what
follows, we need the prior-preconditioned data misfit Hessian, $\ppHmisfit$, defined
in~\eqref{eq:prior-preconditioned-data-misfit-hessian}.  This prior-preconditioned
operator typically exhibits faster spectral decay than
$\Hmisfit$~\cite{Alexanderian_Saibaba_2018}.  This low-rank structure can be used to
construct efficient approximations to the posterior covariance
operator~\cite{Bui-Thanh_Ghattas_Martin_Stadler_2013}.  Here, we use this problem
structure to efficiently approximate the information gain.

Consider a low-rank approximation of $\ppHmisfit$:
\begin{equation}\label{eq:ppHmisfit-lowrank}
  \ppHmisfit \phi
  = \sum_{n=1}^\infty \gamma_n \inp{\phi}{\omega_n} \omega_n
  \approx \sum_{n=1}^r \gamma_n \inp{\phi}{\omega_n} \omega_n,
\end{equation}
where $\phi$ is any function in $\mc{M}$ and $r \leq \Nobs$ is some appropriately chosen
constant such that $(\gamma_n, \omega_n)$ are the dominant eigenpairs of $\ppHmisfit$
given by the eigenvalue problem
\begin{equation} \label{eq:ppHmisfit-eigenproblem}
  \ppHmisfit(\omega_i, \phi)
  = \gamma_i \inp{\phi}{\omega_i},
  \quad \text{with} \quad
  \inp{\omega_i}{\omega_i} = 1
  , \qquad \forall \phi \in \mc{M}, i \in \{1,\dots,r\}.
\end{equation}
We assume that the dominant eigenvalues are not repeated, which is typically the case
due to their rapid decay in inverse problems under study. This assumption is sufficient
to ensure the regularity of the eigenvalue problem, see \cite{Lax_2007}. Using the
low-rank approximation~\eqref{eq:ppHmisfit-lowrank}, we approximate the information gain
as follows:
\begin{equation} \label{eq:aprx-kld}
  \infogain
  \approx 
  \frac{1}{2}
  \left[
  \sum_{i=1}^r \log(1 + \gamma_i)
  - \sum_{i=1}^r \frac{\gamma_i}{1 + \gamma_i}
  + \| \mpost - \mprior \|_{\Cprior^{-1}}^2
  \right].
\end{equation}
This will be our surrogate for the purposes of sensitivity analysis.  Note that the
(low-rank) spectral decomposition of the prior-preconditioned data misfit Hessian is
typically computed in the process of Bayesian analysis in linear(ized) inverse problems.
In such cases, without any additional effort, we have at hand a fast method for
evaluating the information gain.  Such low-rank approximations to the information gain
were also considered, in a discretized setting in~\cite{Alexanderian_Saibaba_2018}.
Here, we utilize this approximation, in an infinite-dimensional setting, to derive
efficient adjoint-based expressions for the eigenvalue derivatives using which we can
efficiently compute the derivatives of the information gain with respect to auxiliary
parameters. 

\subsection{Adjoint-Based Eigenvalue Sensitivity}
\label{sec:eigen-sens}
Since the eigenvalues are only involved in the first two terms of \eqref{eq:aprx-kld},
we first focus on this part for sensitivity analysis and consider the function
\begin{equation} \label{eq:aprx-kld-1}
  \infogain^{\gamma} \defeq \frac{1}{2} \sum_{i=1}^r \left[
    \log(1+\gamma_i) - \frac{\gamma_i}{1+\gamma_i}
    \right].
\end{equation}
This has a dependency on $\vtheta$ through the eigenvalues of the prior-preconditioned
data misfit Hessian $  \ppHmisfit$. Hence, the question of computing
$\partial_{\theta_j}\infogain^{\gamma}$ comes down to computing
$\partial_{\theta_j}\gamma_i$, for $j \in \{1, \ldots, r\}$. To take advantage of the
existing adjoint-based expressions for the data-misfit Hessian, we transform the
eigenvalue problem \eqref{eq:ppHmisfit-eigenproblem} to the generalized eigenvalue
problem
\begin{equation} \label{eq:Hmisfit-eigenproblem}
  \Hmisfit(\psi_i, \phi)
  = \gamma_i \inp{\phi}{\psi_i}_{\Cprior^{-1}}
  \quad \text{and} \quad
  \inp{\psi_i}{\psi_i}_{\Cprior^{-1}} = 1
  , \qquad \forall \phi \in \mc{M}, \forall i \in \{1,\dots,r\},
\end{equation}
where $\psi_i = \Cprior^{1/2}\omega_i$.  Since evaluating \eqref{eq:aprx-kld-1} only
requires eigenvalues which remain unchanged in \eqref{eq:Hmisfit-eigenproblem}, it is
equivalent to use this new eigenvalue problem. Now, leveraging
\eqref{eq:abstract-Hmisfit-system} which implicitly defines the action of $\Hmisfit$, we
can write the full system constraining the evaluation of \eqref{eq:aprx-kld-1} as
\begin{subequations} \label{eq:aprx-kld-1-full}
  \begin{equation} \label{eq:aprx-kld-1-objective}
    \infogain^{\gamma} = \frac{1}{2} \sum_{i=1}^r \left[
      \log(1 + \ms{C}(\hat{p}_i, \psi_i; \vtheta))
      - \frac{
        \ms{C}(\hat{p}_i, \psi_i; \vtheta)
      }{
        1 + \ms{C}(\hat{p}_i, \psi_i; \vtheta)
      }
      \right],
  \end{equation}
  where for $i \in \{1, \dots, r\}$, 
  \begin{equation} \label{eq:aprx-kld-1-constraints}
    \begin{cases}
      \inp{\mcb{Q}\tilde{u}}{\mcb{Q}\hat{u}_i}_{\Cnoise^{-1}}
      + \ms{A}(\hat{p}_i,\tilde{u};\vtheta)
      = 0,                                     & \forall \tilde{u} \in \mc{V}, \\
      \ms{A}(\tilde{p},\hat{u}_i;\vtheta)
      + \ms{C}(\tilde{p},\psi_i;\vtheta)
      = 0,                                     & \forall \tilde{p} \in \mc{V},  \\
      \inp{\psi_i}{\psi_i}_{\Cprior^{-1}} = 1, & \forall \phi \in \mc{M}.
    \end{cases}
  \end{equation}
\end{subequations}
To differentiate~\eqref{eq:aprx-kld-1-full}, we follow a formal Lagrange approach.  To
that end, we construct the meta-Lagrangian
\begin{equation}
  \begin{aligned}
     &
    \mc{L}^\gamma \left(
    \{\psi_i\}_{i=1}^r,
    \{\hat{u}_i\}_{i=1}^r, \{\hat{p}_i\}_{i=1}^r,
    \{\gamma_i^*\}_{i=1}^r,
    \{\hat{u}_i^*\}_{i=1}^r, \{\hat{p}_i^*\}_{i=1}^r
    ; \vtheta
    \right)                                    \\
     & \quad =
    \frac{1}{2} \sum_{i=1}^r \left[
      \log(1 + \ms{C}(\hat{p}_i, \psi_i; \vtheta))
      - \frac{
        \ms{C}(\hat{p}_i, \psi_i; \vtheta)
      }{
        1 + \ms{C}(\hat{p}_i, \psi_i; \vtheta)
      }
    \right]                                    \\
     & \qquad + \sum_{i=1}^r \left[
      \inp{\mcb{Q}\hat{u}_i^*}{\mcb{Q}\hat{u}_i}_{\Cnoise^{-1}}
      + \ms{A}(\hat{p}_i, \hat{u}_i^*; \vtheta)
    \right]                                    \\
     & \qquad + \sum_{i=1}^r \left[
      \ms{A}(\hat{p}_i^*, \hat{u}_i; \vtheta)
      + \ms{C}(\hat{p}_i^*, \psi_i; \vtheta)
    \right]                                    \\
     & \qquad + \sum_{i=1}^r \gamma^*_i \left[
    \inp{\psi_i}{\psi_i}_{\Cprior^{-1}} - 1
    \right].
  \end{aligned}
\end{equation}
Setting the variation with respect to $\hat{p}_i$ in the direction $\tilde{p}$ to zero,
we find:
\begin{equation} \label{eq:eigenvalue-incremental-state}
  \mc{L}^\gamma_{\hat{p}_i}(\tilde{p})
  =
  \frac{
    \ms{C}(\hat{p}_i, \psi_i; \vtheta)
  }{
    2(1 + \ms{C}(\hat{p}_i, \psi_i; \vtheta))^2
  }
  \ms{C}(\tilde{p}, \psi_i; \vtheta)
  + \ms{A}(\tilde{p}, \hat{u}_i^*; \vtheta)
  = 0
  , \qquad \forall \tilde{p} \in \mc{V}.
\end{equation}
Likewise, considering the variation with respect to $\hat{u}_i$ in direction
$\tilde{u}$, we obtain
\begin{equation}\label{eq:eigenvalue-incremental-adjoint}
  \mc{L}^\gamma_{\hat{u}_i}(\tilde{u})
  = \ms{A}(\hat{p}_i^*, \tilde{u}; \vtheta)
  + \inp{\mcb{Q} \tilde{u}}{\mcb{Q} \hat{u}_i^*}_{\Cnoise^{-1}}
  = 0, \qquad \forall \tilde{u} \in \mc{V}.
\end{equation}
Note that \eqref{eq:eigenvalue-incremental-state} and 
\eqref{eq:eigenvalue-incremental-adjoint} are merely rescaled versions of the
incremental state and adjoint equations in~\eqref{eq:aprx-kld-1-constraints}.  In fact, 
\begin{equation} \label{eq:eigenvalue-lagrange-multipliers}
  \hat{u}_i^* = \frac{\gamma_i}{2(1+\gamma_i)^2}\hat{u}_i
  \quad \text{ and } \quad
  \hat{p}_i^* = \frac{\gamma_i}{2(1+\gamma_i)^2}\hat{p}_i
  , \qquad \forall i = 1, \dots, r.
\end{equation}
Thus, no additional PDE solves are needed for $\{ (\hat{u}_i^*, \hat{p}_i^*)\}_{i=1}^r$.
Finally, the derivative of
$\eqref{eq:aprx-kld-1}$ with respect to $\theta_j$ is given by:
\begin{equation} \label{eq:eigenvalue-sensitivity}
  \mc{L}^\gamma_{\theta_j}
  = \sum_{i=1}^r \left[
    \frac{
      \ms{C}(\hat{p}_i, \psi_i; \vtheta)
    }{
      2(1 + \ms{C}(\hat{p}_i, \psi_i; \vtheta))^2
    } \ms{C}_{\theta_j}(\hat{p}_i, \psi_i;\vtheta)
    + \ms{A}_{\theta_j}(\hat{p}_i^*, \hat{u}_i ;\vtheta)
    + \ms{C}_{\theta_j}(\hat{p}_i^*, \psi_i;\vtheta)
    + \ms{A}_{\theta_j}(\hat{p}_i, \hat{u}_i^* ;\vtheta)
    \right].
\end{equation}
Note that all auxiliary parameters are contained in the forms $\ms{A}$ and $\ms{C}$.
Hence, when taking variations in their direction the Lagrange multipliers $\gamma_j^*$,
$j \in \{1, \ldots, r\}$ are not needed.  However, it can be shown that:
\begin{equation}
  \gamma_i^* = \frac{\gamma_i}{2(1+\gamma_i)^2} \gamma_i ,
  \qquad \forall i = 1, \dots, r.
\end{equation}

\subsection{Post-Optimal Sensitivity Analysis}
\label{sec:HDSA}
Here, we consider the third term of \eqref{eq:aprx-kld}. To differentiate this term with
respect to $\theta_j$, $j \in \{1, \ldots, n_\theta\}$, it is sufficient to compute
$\pdv{}{\theta_j}\mpost$ for each $j$. This is possible through a process known as
post-optimal sensitivity analysis or hyper-differential sensitivity analysis (HDSA).  In
particular, previous work
\cite{Sunseri_Hart_Waanders_Alexanderian_2020,Sunseri_Alexanderian_Hart_Waanders_2022}
demonstrates that we are able to compute the derivative of the MAP point $\mpost$ with
respect to any auxiliary parameter $\vtheta$ via the relation
\[
  \pdv{}{\theta_j} \mpost(\vtheta)
  = -\mc{H}^{-1}(\vtheta)\mc{B}_j(\vtheta).
\]
Here, $\mc{B}_j: \Theta_j \to \mc{E}$ is the operator describing the Fr\'echet
derivative of the gradient $\mc{G}$, defined in~\eqref{eq:abstract-gradient-system},
with respect to $\theta_j$. $\mc{H}^{-1}$ is the inverse Hessian operator, and its
action is discussed later in \cref{sec:alg}. This derivative is then evaluated at
$\mpost$. Differentiating the third term of the information gain \eqref{eq:aprx-kld},
we find:
\begin{equation} \label{eq:kld-hdsa}
  \pdv{}{\theta_j} \| \mpost - \mprior \|_{\Cprior^{-1}}^2
  = 2 \inp{ \mpost - \mprior} {\pdv{}{\theta_j} \mpost }_{\Cprior^{-1}}
  = -2 \inp{\mpost - \mprior} {\mc{H}^{-1} \mc{B}_j}_{\Cprior^{-1}}.
\end{equation}
Next, we derive the adjoint-based expression for $\mc{B}_j$.  This is done again via a
Lagrange multiplier approach.  Namely, we differentiate through $\mc{G}$, by setting up
a meta-Lagrangian where we enforce the  state and adjoint equations via Lagrange
multipliers.  In fact, we consider the same meta-Lagrangian constructed for the Hessian
computation~\eqref{eq:lagrangian-hessian}, which we restate for clarity:
\begin{equation}\label{eq:lagrangian-hdsa}
\begin{aligned}
  \mc{L}^{\rm HDSA}(
  u, p, m,
  \hat{u}, \hat{p}, \hat{m};
  \vtheta
  )
   & =
  \inp{m-\mprior}{\hat{m}}_{\Cprior^{-1}}
  + \ms{C}(p,\hat{m};\vtheta) \\
   & \quad
  + \inp{\mcb{Q}u-\obs}{\mcb{Q}\hat{u}}_{\Cnoise^{-1}}
  + \ms{A}(p, \hat{u};\vtheta) \\
   & \quad
  + \ms{A}(\hat{p},u;\vtheta)
  + \ms{C}(\hat{p},m;\vtheta)
  + \ms{D}(\hat{p}; \vtheta).
\end{aligned}
\end{equation}
Letting the variations with respect to $u$ and $p$ vanish, we recover the incremental
state and adjoint equations found previously in
\eqref{eq:abstract-incremental-state}-\eqref{eq:abstract-incremental-state}, with
solution $(\hat{u}, \hat{p})$.  Subsequently, differentiating through with respect to
the auxiliary model parameters, we obtain
\begin{equation} \label{eq:post-optimal-B_j}
  \mc{L}^{\text{HDSA}}_{\theta_j}
  = \mc{B}_j(\vtheta)(\hat{m})
  \defeq \ms{C}_{\theta_j}(p, \hat{m}; \vtheta)
  + \ms{A}_{\theta_j}(p, \hat{u}; \vtheta)
  + \ms{C}_{\theta_j}(\hat{p}, m; \vtheta)
  + \ms{A}_{\theta_j}(\hat{p}, u; \vtheta)
  + \ms{D}_{\theta_j}(\hat{p}; \vtheta),
\end{equation}
where $\mc{A}_{\theta_j}, \mc{C}_{\theta_j}$, and $\mc{D}_{\theta_j}$ are the derivative
of the given weak form with respect to $\theta_j$.  In practice, the assembly of
$\mc{B}_j$ requires two additional incremental state and adjoint variable solves.  The
details of this construction are explored in
\cite{Sunseri_Alexanderian_Hart_Waanders_2022}.

\subsection{The overall algorithm and computational costs}
\label{sec:alg}

The developments in the previous subsections provide the building blocks for computing
the sensitivity of the information gain with respect to auxiliary parameters.
\Cref{alg:local_sensitivity} summarizes our overall approach. In practice, the
computational burden for this algorithm is contained in a few portions of the procedure.
Below, as well as in table \Cref{table:costs}, we summarize the computational cost in
the number of PDE solves required for each step. These solves typically dominate the
computation, and vary heavily depending on the discretization, method employed, and
individual characteristics of the equation.

In infinite dimensions the prior covariance operator is typically defined as the inverse
of a differential operator.  Thus, computing the action of the prior covariance operator
requires PDE solves. We do not include this in our computational cost analysis, because
the prior covariance typically involves a constant coefficient scalar elliptic operator
which can be pre-factorized, while the state equation typically is a vector system,
often with spatially varying coefficients, and potentially time-dependent. Thus, solving
the state equation (and its adjoint and incremental variants) dominates the cost of
$\Cprior$ applies. 
\begin{algorithm}
  \caption{
    Computing local sensitivities of the information gain with respect to $\theta_j$
  }\label{alg:local_sensitivity}
  \DontPrintSemicolon
  \LinesNumbered
  \SetNoFillComment
  \KwData{
  $\{(\gamma_i, \omega_i)\}_{i=1}^r$ and MAP point $\mpost$
  }
  \KwResult{
    $S \defeq \pdv{}{\theta_j}\infogain^{\gamma}$
  }
  %\BlankLine
  $S \gets 0$\;
  Solve the state and adjoint equations for $(u, p)$:
  \[
    \begin{aligned}
    \ms{A}(\tilde{p}, u; \vtheta)
    + \ms{C}(\tilde{p}, \mpost; \vtheta)
    + \ms{D}(\tilde{p}; \vtheta)
    &= 0, \quad \forall \tilde{p} \in \mc{V}\\
    \inp{\mcb{Q} u - \obs}{\mcb{Q} \tilde{u}}_{\Cnoise^{-1}}
    + \ms{A}(\hat{p}, \tilde{u}; \vtheta)
    &= 0, \quad \forall \tilde{u} \in \mc{V}
    \end{aligned}
  \]\; \vspace{-0.5cm}
  Solve the incremental state and incremental adjoint equations for $(\hat{u},
  \hat{p})$:
  \[
    \begin{aligned}
      \ms{A}(\tilde{p}, \hat{u}; \vtheta)
      + \ms{C}(\tilde{p}, \mpost; \vtheta)
      &= 0, \quad \forall \tilde{p} \in \mc{V} \\
      \inp{\mcb{Q} \tilde{u}}{\mcb{Q} \hat{u}}_{\Cnoise^{-1}}
      + \ms{A}(\hat{p}, \tilde{u}; \vtheta)
      &= 0, \quad \forall \tilde{u} \in \mc{V}
    \end{aligned}
  \]\;\vspace{-0.5cm}
  \For{$i \in \{1, \dots, r\}$}{
    \tcc{Eigenvalue Sensitivity}
    $\psi_i \gets \Cprior^{1/2} \vec{\omega}_i$
    \tcp*{Convert to eigenfunction of $\Hmisfit$}
    Solve the incremental state and incremental adjoint (in direction $\psi_i$) to
    obtain $(\hat{u}_i, \hat{p}_i)$:
    \[
      \begin{aligned}
        \ms{A}(\tilde{p}, \hat{u}; \vtheta)
        + \ms{C}(\tilde{p}, \psi_i; \vtheta)
        &= 0, \quad \forall \tilde{p} \in \mc{V} \\
        \inp{\mcb{Q} \tilde{u}}{\mcb{Q} \hat{u}}_{\Cnoise^{-1}}
        + \ms{A}(\hat{p}, \tilde{u}; \vtheta)
        &= 0, \quad \forall \tilde{u} \in \mc{V}
      \end{aligned}
    \]\;\vspace{-0.5cm}
    Construct Lagrange multipliers $(\hat{u}_i^*, \hat{p}_i^*) =
    \frac{\gamma_i}{2(1+\gamma_i)^2}(\hat{u}_i, \hat{p}_i)$\;
    Let
    \(
      S \gets S
      + \frac{\gamma_i}{2(1+\gamma_i)^2}
      \ms{C}_{\theta_j}(\hat{p}_i, \psi_i;\vtheta)
      + \ms{A}_{\theta_j}(\hat{p}_i^*, \hat{u}_i ;\vtheta)
      + \ms{C}_{\theta_j}(\hat{p}_i^*, \psi_i;\vtheta)
      + \ms{A}_{\theta_j}(\hat{p}_i, \hat{u}_i^* ;\vtheta)
    \)\;
    \tcc{Post-Optimal Sensitivity}
    Build
    \(
      \mc{B}_j
      = \ms{C}_{\theta_j}(p, \hat{m}; \vtheta)
      + \ms{A}_{\theta_j}(p, \hat{u}; \vtheta)
      + \ms{C}_{\theta_j}(\hat{p}, m; \vtheta)
      + \ms{A}_{\theta_j}(\hat{p}, u; \vtheta)
      + \ms{D}_{\theta_j}(\hat{p}; \vtheta)
    \)\;
    Let $S \gets S - \inp{m - \mprior}{\mc{H}^{-1}\mc{B}_j}_{\Cprior^{-1}}$
  }
  \Return{S}
\end{algorithm}
The Hessian is a central point of consideration, as nearly all presented methods
mentioned revolve around its action. In light of its low-rank structure, we begin by
considering the cost of computing the spectral decomposition of the prior-preconditioned
data-misfit Hessian. The Hessian action itself requires the solution of the incremental
state and adjoint PDEs; cf.~\eqref{eq:abstract-Hmisfit-system}.  Hence, if we employ a
Krylov iterative method, requiring $\mc{O}(r)$ matrix vector applications to compute a
rank $r$ spectral decomposition, $\mc{O}(2 \times r)$ PDE solves are needed.  Note that
this rank $r$ is bounded by the number of data observations, which is independent of
discretization dimension.

Likewise, the application of the inverse Hessian is central to computing the MAP point
as well as the post-optimal sensitivity analysis step for each parameter. Note, that for
the linear Gaussian inference setting, the inverse Hessian is the posterior covariance
and has a closed form expression given by
\begin{equation}
  \mc{H}^{-1}
  = (\Hmisfit + \Cprior^{-1})^{-1}
  = \Cprior^{1/2}(\ppHmisfit + \mc{I})^{-1}\Cprior^{1/2}.
\end{equation}
It has been demonstrated, see e.g., \cite{Alexanderian_Gloor_Ghattas_2016}, that
\begin{equation}
  (\ppHmisfit + \mc{I})^{-1}z
  = z - \sum_{j=1}^\infty \frac{\gamma_j}{1+\gamma_j} \inp{z}{\omega_j}_{\mc{M}} \omega_j
  \approx z - \sum_{j=1}^r \frac{\gamma_j}{1+\gamma_j} \inp{z}{\omega_j}_{\mc{M}} \omega_j
  , \quad \forall z \in \mc{M},
\end{equation}
where $(\gamma_j, \omega_j)_{j=1}^\infty$ are the eigenpairs of $\ppHmisfit$.
Therefore, after the spectral decomposition of $\ppHmisfit$, the approximate application
of the inverse Hessian requires no additional PDE solves. This implies, for example,
that the MAP point estimation requires no additional PDE solves.

Considering the local sensitivity calculation, naturally the work required will scale
according to $n_\theta$.  However, in terms of PDE solves, this is only the case for the
post-optimal sensitivity analysis.  Indeed, for each of the $n_\theta$ parameters, one
must solve the incremental state and adjoint systems to evaluate
$\eqref{eq:post-optimal-B_j}$; hence $\mc{O}(2 \times n_\theta)$ PDE solves.  However,
for the eigenvalue sensitivities, \eqref{eq:eigenvalue-sensitivity} requires the
solution of $\mc{O}(2 \times r)$ PDEs to obtain the incremental state and adjoint
variables corresponding to each eigenfunction, as given by
\eqref{eq:eigenvalue-incremental-state} and \eqref{eq:eigenvalue-incremental-adjoint}.
These incremental variables are then reused for the sensitivity calculation for each
$\theta_j, j = 1, \dots, n_\theta$. Hence, the cost in number of PDE solves for the
eigenvalue sensitivity calculation is $\mc{O}(2 \times r)$, independent of $n_\theta$

As a brief note on the global sensitivity analysis procedure, equation
\eqref{eq:sobol-bound} would typically be evaluated using sample average style methods.
That is, for sufficiently large $\Ns$:
\begin{align}
  \label{eq:variance-sample-average}
  \Var{\infogain(\vtheta)} 
  &\approx
  \frac{1}{\Ns-1}\sum_{i=1}^{\Ns} \left( 
    \infogain(\vtheta_i) - \frac{1}{\Ns}\sum_{i=1}^{\Ns} \infogain(\vtheta_i)
  \right)^2
  , \\
  \label{eq:dgsm-sample-average}
  \int \left(
    \pdv{}{\theta_i}\infogain(\vtheta)
  \right)^2 \dd{F(\vtheta)}
  &\approx
  \frac{1}{\Ns} \sum_{i=1}^{\Ns} \pdv{}{\theta_i}\infogain(\vtheta_i)^2.
\end{align}
This would require \Cref{alg:local_sensitivity} to be run $\Ns$ times and the inverse
problem itself to be solved $\Ns$ times. Note, these estimators may require a different
number of samples to converge. However, the strength of \Cref{alg:local_sensitivity} is
its ``value \textit{and} gradient" operation, as the sensitivity calculation takes
advantage of low ranking already computed during the Bayesian inversion procedure.
Therefore, it is convenient to reuse the samples between the two estimators.

\begin{table}[!htb]
  \centering
  \begin{tabular}{l|l}
    \textbf{Computation}           & \textbf{PDE Solves Required}              \\ \hline
    Low-Rank $\ppHmisfit$ Approximation & $\mc{O}(2 \times r)$ \\ \hline
    MAP Point Estimation          & $\mc{O}(1)$ \\ \hline
    Eigenvalue Sensitivities       & $\mc{O}(2\times r)$ \\ \hline
    Post-Optimal Sensitivity       & $\mc{O}(2 \times n_\theta)$ \\
  \end{tabular}
  \caption{
    Computational cost summary.
  }\label{table:costs}
\end{table}

\section{Computational Examples}\label{sec:numerics}
We consider two computational examples.  The first one, discussed
in~\cref{sec:elliptic-toy}, involves source inversion in a (scalar) linear elliptic PDE.
Here, we illustrate the use of information gain and expected information gain as HDSA
QoIs.  The second inverse problem, presented in~\cref{sec:cr-model}, is governed by the
equations of linear elasticity in three dimensions. We fully elaborate our proposed
framework for this example, which is motivated by a geophysics inverse problem that aims
to infer the slip field along a fault after an earthquake. Both of the following
numerical examples were built using the PDE-constrained inverse problem capabilities of
the open-source software package \texttt{hIPPYlib}~\cite{VillaPetraGhattas21}.

\subsection{Source Inversion in an Elliptic PDE} 
\label{sec:elliptic-toy}

On the domain $\Omega := (0,1)^2$, we consider the elliptic PDE
\begin{subequations} \label{eq:elliptic-toy}
  \begin{alignat}{2}
    -\Delta u + c u & = m, \quad \text{in } \Omega  ,       \\
    \nabla u \cdot n  & = g, \quad \text{on } \partial\Omega.
  \end{alignat}
\end{subequations}
We aim to infer the volume source term $m$ from point observations of the solution $u$.
We consider the parameters $c$ and $g$ as scalar auxiliary parameters with respect to
which we study the sensitivity of the solution of the Bayesian inverse problem.  In the
present study, the nominal values of these parameters are $c = 1$ and $g = 0.1$.  To
obtain synthetic observation data, we use a true source term $m^{\rm true}$, given by
\begin{equation}\label{eq:elliptic-truem}
  m^{\rm true}(x,y) = 10 e^{-\frac{1}{20}\left[(x-0.5)^2 + (y-0.5)^2\right]}.
\end{equation}
We then solve \eqref{eq:elliptic-toy} with this synthetic source to obtain $u^{\rm
true}$ using a finite element discretization with first-order Lagrangian triangular 
elements on a $32 \times 32$ grid. The observations are obtained by evaluating the
solution at $9$ evenly distributed points in $\Omega$, along with adding additive
Gaussian noise with standard deviation $\sigma = 0.01 \| u \|_\infty$. As prior
distribution, we use a Gaussian prior with zero mean and a squared inverse elliptic
operator $\Cprior = (I - \Delta)^{-2}$ as described in~\cite{VillaPetraGhattas21}.

\begin{figure}[!htb]
  \centering
  \begin{subfigure}[b]{0.495\textwidth}
    \caption{
      Expected information gain $\Einfogain$
    }
    \vspace{5pt}
    \includegraphics[width=\textwidth]{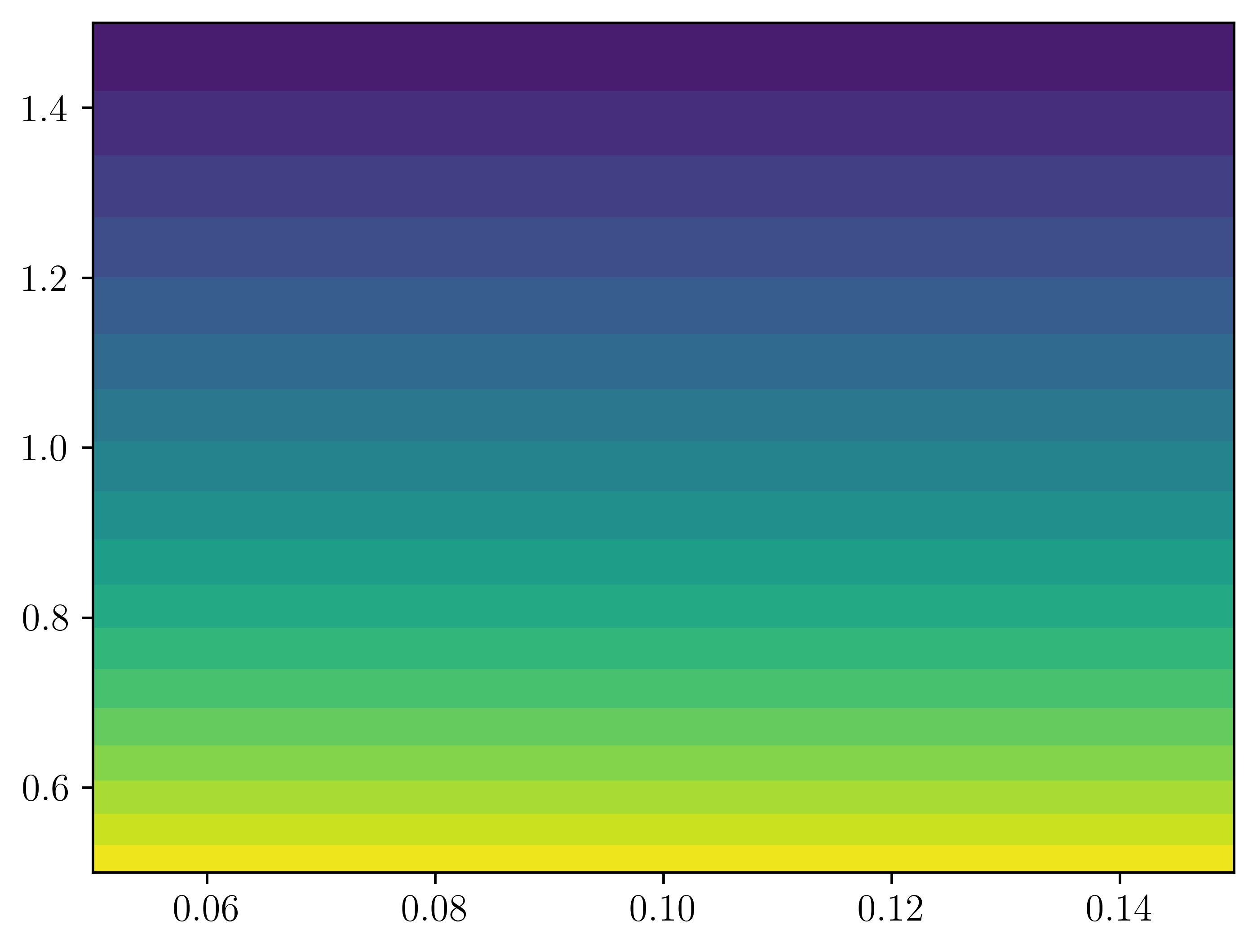}
    \includegraphics[width=\textwidth]{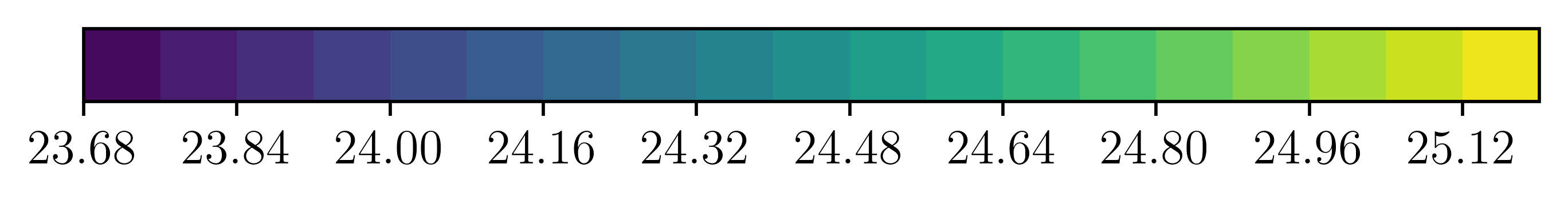}
  \end{subfigure}
  \begin{subfigure}[b]{0.495\textwidth}
    \caption{
      Information Gain $\infogain$
    }
    \vspace{5pt}
    \includegraphics[width=\textwidth]{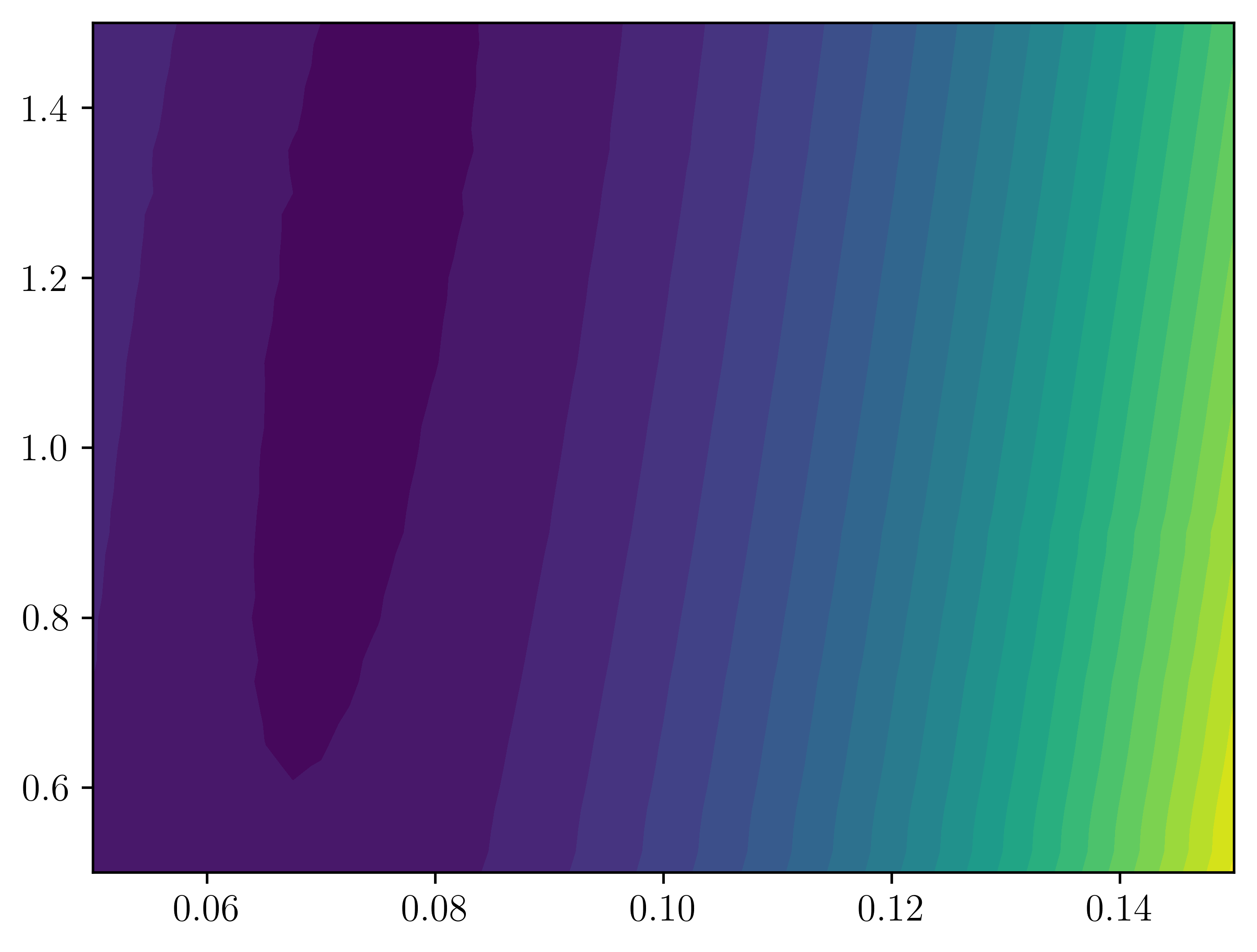}
    \includegraphics[width=\textwidth]{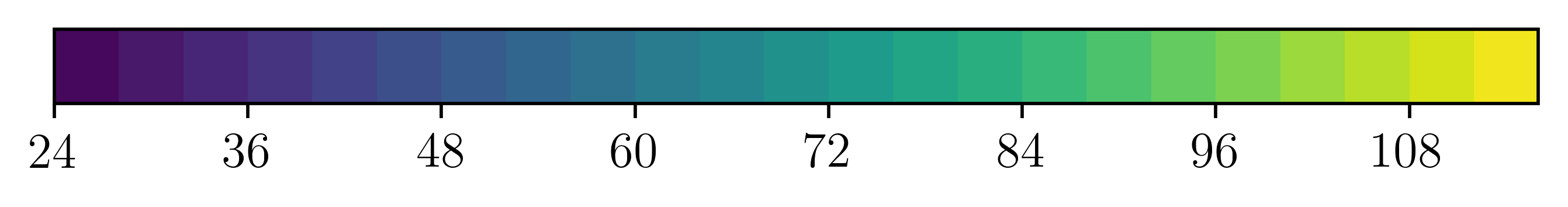}
  \end{subfigure}
  \caption{
    $\infogain$ and $\Einfogain$ as functions of $g$ ($x$-axis) and $c$ ($y$-axis). Note
    the difference in scales and that $\Einfogain$ is independent of $g$.
  }\label{fig:elliptic_infogain}
\end{figure}

With the inverse problem fully specified, we invoke the machinery developed in
Section~\ref{sec:methods} to compute the information gain and its derivatives at the
nominal auxiliary parameter values.  In particular, we use this example to study the
dependence of information gain as well as expected information gain to $g$ and $c$. To
this end, in \Cref{fig:elliptic_infogain}, we report the expected information gain
(left) and information gain (right) as functions of the auxiliary parameters. 
We note that the information gain varies strongly as a function of $g$. This sensitivity
is not unexpected from this model, because $g$ influences flow through the boundary. 
Moreover, since $g$ is a boundary source term that enters the governing equations
linearly, it does not appear in the Hessian \eqref{eq:abstract-hessian-action}.
Therefore, the influence of $g$ on the information gain is solely coming from the third
term in~\eqref{eq:linear-kld}.  In particular, the expected information gain is
independent of $g$ as can be seen in \Cref{fig:elliptic_infogain}.  Moreover, note that
$\infogain$, when considered as a function of $g$, has a local minimum at which the
sensitivity vanishes.  From observing this local sensitivity, one may draw the incorrect
conclusion that $\infogain$ is generally insensitive to $g$. In such a case, a global
sensitivity analysis, as demonstrated in the next section, is a more accurate
description.

\subsection{Fault-Slip Inversion in 3D} \label{sec:cr-model}
In this section, we apply the proposed methods to a three-dimensional model
problem inspired from the field of seismology. Namely, we consider the setting
of linear elasticity to describe the deformation of the earth's crust due to
tectonic activity. An important inverse problem in this context is that of
inferring the slip pattern along a tectonic fault based on observations of the
surface deformation.  The specific problem we study is adapted
from~\cite{McCormack_2018}. Although our model is merely a synthetic example, it
is representative of the type of inverse problem our methods are intended to
analyze.

\subsubsection{Forward Model}
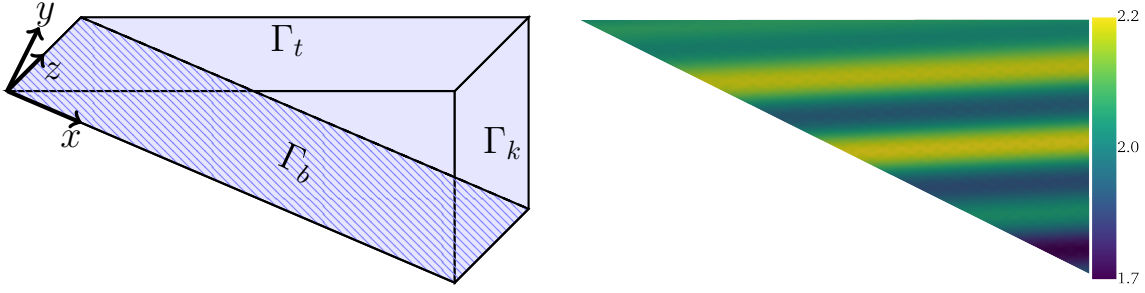
\begin{figure}[!htb]
  \centering
  \begin{tikzpicture}[rotate=0,scale=0.85]
    \pgfmathsetmacro{\cubex}{7}
    \pgfmathsetmacro{\cubey}{3}
    \pgfmathsetmacro{\cubez}{3}
    \draw[thick,black,fill=blue!10] (-\cubex,0,0) -- ++(\cubex,-\cubey,0) -- (0,0,0) -- cycle;
    \draw[thick,black,fill=blue!10] (0,0,0) -- ++(0,0,-\cubez) -- ++(0,-\cubey,0) -- ++(0,0,\cubez) -- cycle;
    \draw[thick,black,fill=blue!10,pattern color=blue!50] (0,0,0) -- ++(-\cubex,0,0) -- ++(0,0,-\cubez) -- ++(\cubex,0,0) -- cycle;
    \draw[dashed,thick,black,fill=blue!10] (-\cubex, 0, -\cubez) -- ++(\cubex,-\cubey,0) -- cycle;
    \draw[thick,pattern=north west lines,pattern color=blue!50] (-\cubex,0,0) -- (-\cubex,0,-\cubez) -- (0,-\cubey,-\cubez) --
                            (0,-\cubey,0) -- cycle;
    \coordinate [label=\Large$\Gamma_t$] (dom) at (-3,0,-1);
    \node[label=below:\rotatebox{-25}{\Large$\Gamma_b$}] at (1,3,9) {};
    \node[label=below:\rotatebox{0}{\Large$\Gamma_k$}] at (0.75,-0.25,0) {};
    %
    % add axes
    %
    \draw[ultra thick,black,->] node at (-6,-.75,0){\Large$x$} (-\cubex,0,0) -- ++(\cubex/6,-\cubey/6,0);
    \draw[ultra thick,black,->] node at (-\cubex+.1,-0.3,-\cubez/2-.1){\Large$z$} (-\cubex,0,0) -- (-\cubex,0,-\cubez/2);
    \draw[ultra thick,black,->] node at (-\cubex+.6,1.2,0){\Large$y$} (-\cubex,0,0) -- (-\cubex+.5,1,0);
  \end{tikzpicture}\hspace{3ex}
  \includegraphics[width=0.45\textwidth]{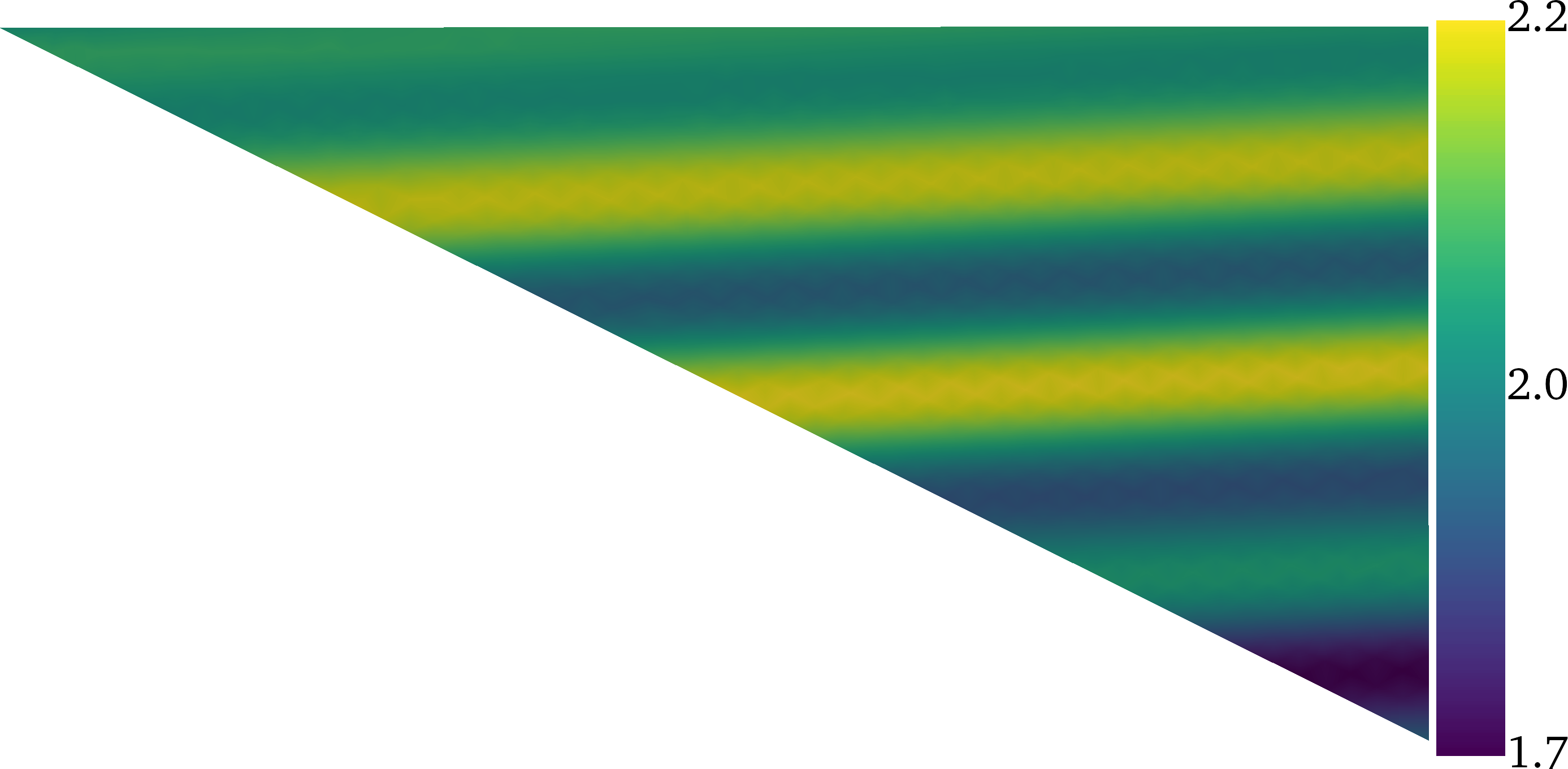}
  \caption{
    (Left) Visualization of $\Omega$. (Right) Realization of $\lambda$ with six KLE
    modes. 
  }\label{fig:lambda}
\end{figure}

Consider a triangular prism domain $\Omega\subset \mathbb R^3$ with
boundary $\Gamma = \Gamma_b\cup \Gamma_k \cup \Gamma_t\cup \Gamma_s$, where
$\Gamma_b, \Gamma_k, \Gamma_t$, and $\Gamma_s$ are the bottom, back, top and the
side surfaces, respectively; see \Cref{fig:lambda} for a visualization of
the domain.  We would like to infer the slip along $\Gamma_b$, denoted by $\vec{m}$, based
on point observations of the displacement on the top surface $\Gamma_t$. We assume a
linear elasticity model for the displacement $\vec{u} = [u_1 \; u_2 \; u_3]^\top \in
H^1(\Omega)^3$:
\begin{equation}
  -\nabla \cdot \vec{\sigma}(\vec{u}) = \vec{0}
  ,\qquad \text{ in } \Omega,
\end{equation}
where $\vec{\sigma}(\vec{u}) := \mathbb{C} \vec{\e}(\vec{u}) = 2\mu \vec{\e(\vec{u})}
+ \lambda \text{tr}(\vec{\e}(\vec{u})) \mat{I}$ is the fourth-order linear elasticity
stress tensor and $\vec\e (\vec{u}) := \frac 12 \left(\nabla \vec{u} + (\nabla
\vec{u})^T \right)$ is the strain tensor. The stress tensor involves the L\'{a}me
parameters $\{\mu, \lambda\}$ describing material properties. Specifying the
boundary
conditions, we have:
\begin{subequations}\label{eq:elasticity_strong}
  \begin{alignat}{2}
    \label{eq:elasticity_strong1}
    -\nabla \cdot [
      \mu(\nabla \vec{u} + (\nabla \vec{u})^T)
      + \lambda\nabla \cdot \vec{u} \mat{I}
    ]                                  & = \vec{0} &  & \text{ in } \Omega,   \\
    \label{eq:elasticity_strong2}
    \vec{\sigma}(\vec{u})\vec{n}       & = 0       &  & \text{ on } \Gamma_t, \\
    \label{eq:elasticity_strong3}
    \vec{u} + \nu_k \vec{\sigma(u)n} & = \vec{0} &  & \text{ on } \Gamma_k, \\
    \label{eq:elasticity_strong4}
    \vec{u} + \nu_s \vec{\sigma(u)n} & = \vec{h} &  & \text{ on } \Gamma_s, \\
    \label{eq:elasticity_strong5}
    \vec{u} \cdot \vec{n}              & = 0       &  & \text{ on } \Gamma_b, \\
    \label{eq:elasticity_strong6}
    \delta \mat{T} (\vec{\sigma}(\vec{u})\vec{n}) + \mat{T}\vec{u}
                                       & = \vec{m} &  & \text{ on } \Gamma_b.
  \end{alignat}
\end{subequations} 
Here, $\mat{T}$ is the operator that extracts the tangential components of a vector,
i.e., $ \mat{T}\vec{u} \defeq (\mat I - \vec{n} \otimes \vec{n})\vec{u} = \vec{u} -
(\vec{n}^T \vec{u})\vec{n}$, where $\vec{n}$ is the unit normal of the domain.  The
Robin-type condition \eqref{eq:elasticity_strong6} can be understood as regularized
Dirichlet condition with parameter $0 < \delta \ll 1$. 

To model the heterogeneous properties of the earth, we use spatially varying L\'{a}me
parameters. We assume these parameters to have variations only in vertical direction
modeling a layered earth structure. Specifically, we define the L\'{a}me parameters
using the expansions
\begin{equation} \label{eq:kl-expansion}
  \lambda(y)
  = \overline{\lambda}
  + \sum_{i=1}^6 e_i(y)\lambda_i^{\rm KLE}
  \qquad \text{and} \qquad
  \mu(y)
  = \overline{\mu}
  + \sum_{i=1}^6 e_i(y)\mu_i^{\rm KLE}.
\end{equation}

This representation is motivated by the truncated Karhunen--Lo\`eve expansion (KLE),
hence the choice of superscript $\mathrm{KLE}$ in the above representation.  In these
equations, $e_i$'s are suitably scaled KLE basis functions, computed as described
in~\cite{LeMatre_Knio_2010}.  For the purposes of the present study, the coefficients
$\{\lambda_i^{\rm KLE}\}$ and $\{\mu_i^{\rm KLE}\}$ are (fixed) random draws from a
standard normal distribution. To illustrate, we show the corresponding realization of 
the $\lambda$ field in~\Cref{fig:lambda}~(right).  As for the remaining auxiliary
parameters, we assume the nominal values of $\overline{\lambda} = 2.0, \overline{\mu} =
2.50, \nu_k = 10^{-4}$, and $\nu_s = 10^{-4}$.

\subsubsection{Inverse Problem Setup}

Now we return to the inverse problem under study. We seek to infer the fault
slip $\vec{m}$ and compute the sensitivity of the information gain with respect
to the auxiliary parameters $\nu_k, \nu_s, \overline{\mu}, \{\mu_i^{\rm KLE}\},
\overline{\lambda}$, and $\{\lambda_i^{\rm KLE}\}$.  To cast the problem into the
abstract setup developed in Section~\ref{sec:methods}, we consider the weak
form of the governing PDE. 
Defining the function space $\vec{V} \defeq \{ \vec{u} \in
H^1(\Omega)^3 : \vec{u} \cdot \vec{n} = 0 \text{ on } \Gamma_b \}$, the weak form
of \eqref{eq:elasticity_strong} is as follows: find $\vec{u} \in \mat{V}$ such that for every
$\vec{p} \in \vec{V}$
\begin{equation*} %\label{eq:elast_weak}
  \int_{\Gamma_k} \nu^{-1}_k\vec{u} \cdot \vec{p} \dd{s}
  + \int_{\Gamma_s} \nu^{-1}_s\vec{u} \cdot \vec{p} \dd{s}
  + \int_{\Gamma_b} \delta^{-1}\mat{T}\vec{u} \cdot \vec{p} \dd{s}
  + \int_\Omega \e(\vec{u}) : \mb{C}[\e(\vec{p})] \dd{x}
  - \int_{\Gamma_b} \delta^{-1} \vec{m} \cdot \vec{p} \dd{s}
  - \int_{\Gamma_s} \nu_s^{-1}\vec{h} \cdot \vec{p} \dd{s}
  = 0.
\end{equation*}
Note that this fits into the abstract weak PDE formulation described in Section
\ref{sec:math-prelim} (see~\eqref{eq:abstract-pde}), and facilitates deploying the
approach developed in Section~\ref{sec:methods}.

\begin{figure}[ht]
  \centering
  \includegraphics[width=0.49\textwidth]{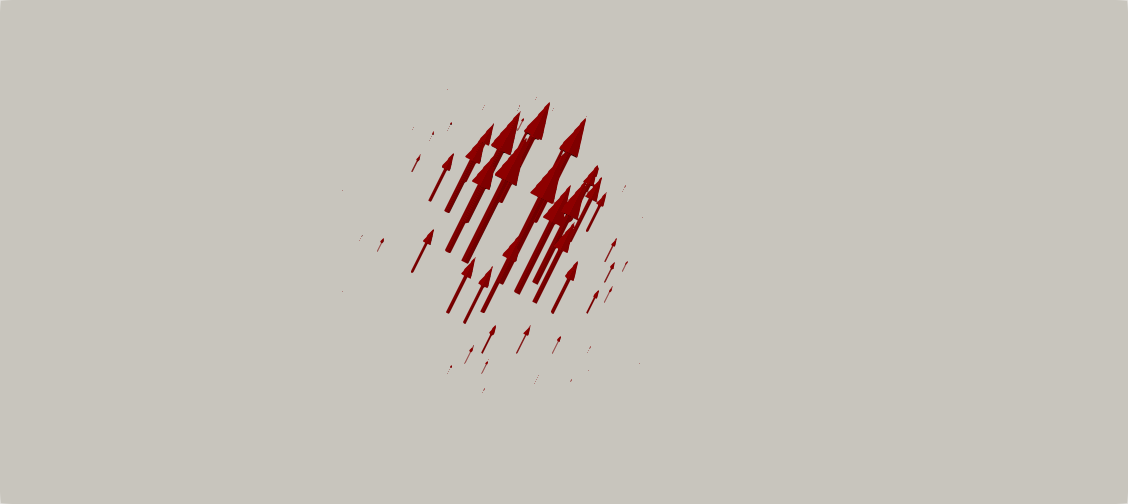}
  \includegraphics[width=0.49\textwidth]{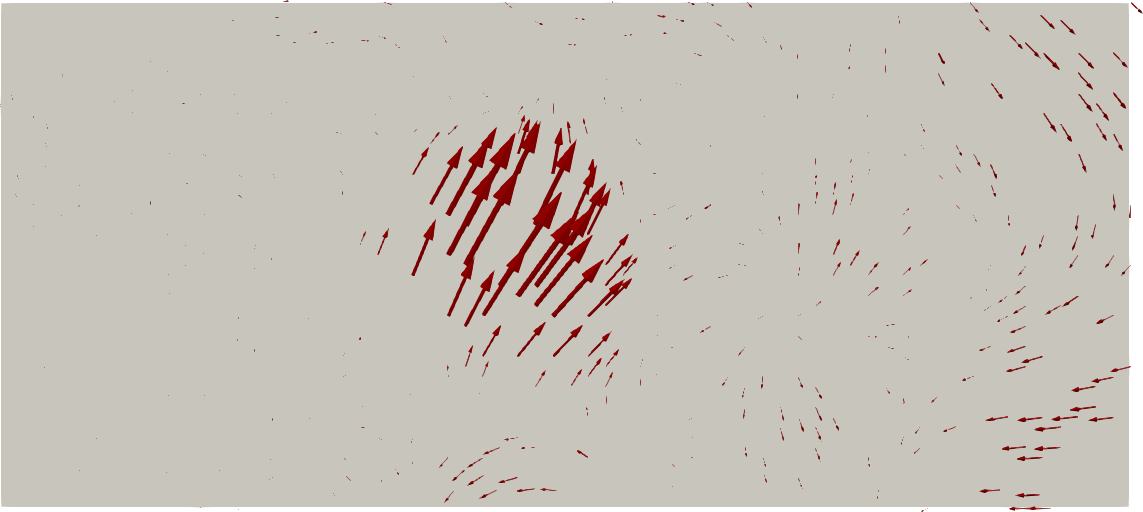}
  \caption{True slip field $\vec m$ given by \eqref{eq:true-slip} (left) and
    the slip MAP estimate $\vec{m}$ (right).} 
  \label{fig:u0}
\end{figure}
\begin{figure}[ht]
  \centering
  \includegraphics[width=0.49\textwidth]{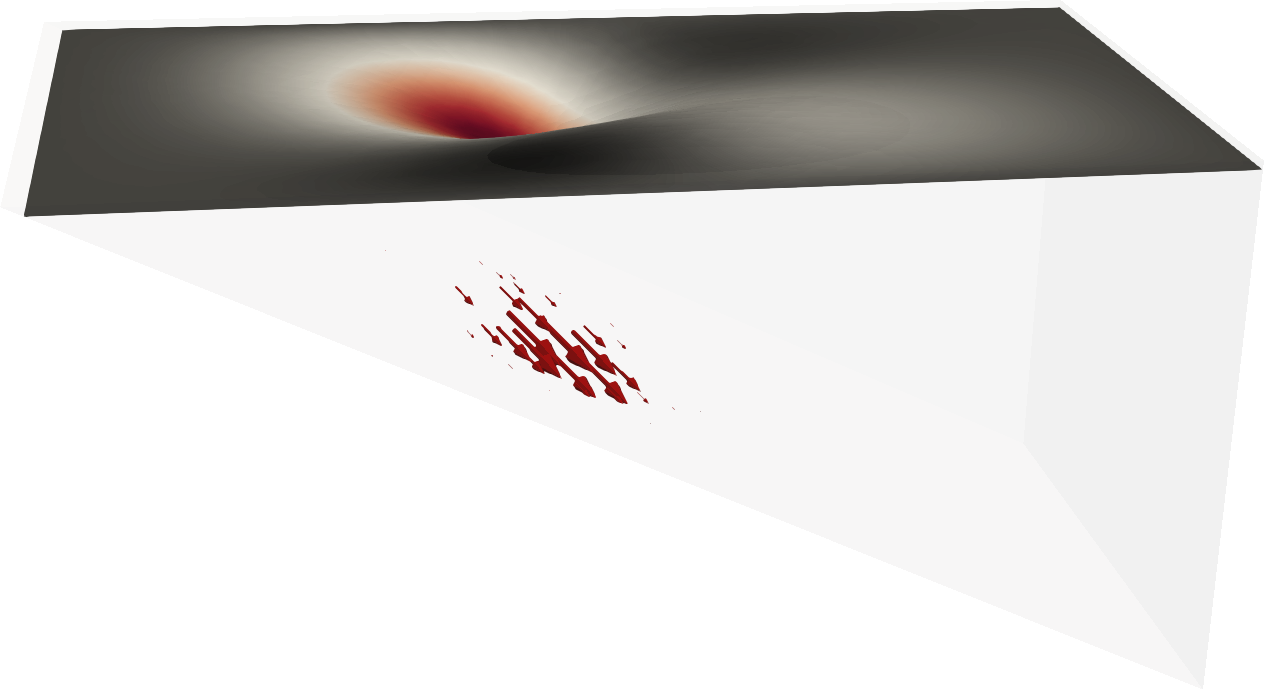}
  \includegraphics[width=0.49\textwidth]{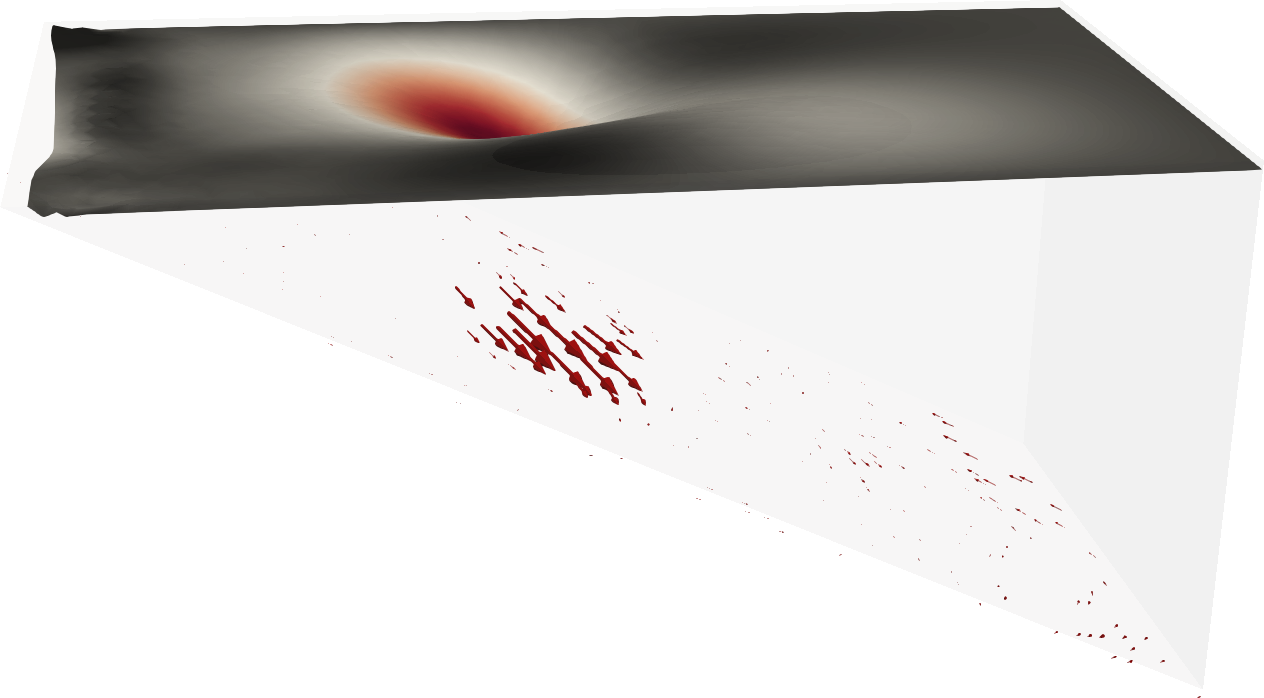}
  \caption{
    Fault slip (red arrows) and corresponding deformation on the top surface $\Gamma_t$.
    Shown on the left is the result for the truth fault slip used to generate synthetic
    data, and on the right the MAP fault slip.
  } \label{fig:state}
\end{figure}

For the synthetic data generation, we assume that the
``ground-truth'' fault slip has the form
\begin{equation}\label{eq:true-slip}
  \vec{m}(x,y,z) = \begin{pmatrix}
    \exp\left[(x+1.25)^2 + (z - 0.5)^2\right] \\
    0                                             \\
    2 \exp\left[(x+1.25)^2 + (z - 0.5)^2\right] \\
  \end{pmatrix}.
\end{equation}
See \Cref{fig:u0} (left) for a visualization of this ground-truth parameter.  We solve
our forward model \eqref{eq:elasticity_strong} with this slip field and our nominal set
of parameters to compute the ``truth'' displacement field $\vec{u}^{\rm true}$, which is
shown on the left in \Cref{fig:state}.  We then take point measurements on $\Gamma_t$,
the top surface of the domain at $\Nobs = 64$ measurement locations, spaced uniformly on
$\Gamma_t$. We add independent additive noise with the standard deviation $\sigma =
10^{-3}$ to this measurement data to obtain synthetic observations $\obs$.

Next, we discuss the prior measure. We use a Gaussian measure with mean $\mprior =
\vec{0}$, and a squared inverse elliptic operator as the covariance. Specifically,
$\Cprior = \mc{S}^{-2}$, where  $\mc{S}$ is the solution operator mapping $\vec{s}$ to
$\vec m$ for the following PDE given in weak form as:
\begin{equation} \label{eq:elasticity-prior}
  \int_{\Gamma_b}
  \gamma_\text{pr} \nabla \vec{m} \cdot \nabla \vec{p} + \delta_\text{pr} \vec{m} \cdot \vec{p}
  \dd{x} 
  + \int_{\partial \Gamma_b} \beta_\text{pr} \vec{m} \cdot \vec{p} \dd{s}
  = \int_{\Gamma_b} \vec{s} \cdot \vec{p} \dd{x},
  \qquad \forall \vec{p} \in H^1(\Gamma_b)^2.
\end{equation}
In this formulation, $\gamma_\text{pr}$ and $\delta_\text{pr}$ govern the correlation
length and pointwise variance of the prior samples, and $\beta_\text{pr}$ is a parameter
stemming from an additional Robin boundary condition that is used to minimize boundary
artifacts. This formulation of the prior is adapted from~\cite{VillaPetraGhattas21}, and
we choose the parameters $(\gamma_\text{pr}, \delta_\text{pr}) = (0.01, 0.8)$ and
$\beta_\text{pr} = \sqrt{\gamma_\text{pr} \delta_\text{pr}}/1.42$. 

\subsubsection{Numerical Results} \label{sec:numerical-results}
We now solve the Bayesian inverse problem. In \Cref{fig:u0}, we show the MAP point,
which differs from the true slip field as a consequence of the sparse measurements, the
noise in the measurement data, and the prior assumptions. However, the solution of the
Bayesian inverse problem is not the main focus here. Instead, we seek to understand the
sensitivity of the information gain with respect to the auxiliary model parameters,
using \Cref{alg:local_sensitivity}. However, the computation of the information
gain and its sensitivity requires specification of the hyperparameters of the low-rank
approximation to the prior-preconditioned data-misfit Hessian. As a natural
upper bound for the rank of the Hessian is the number of observations, one could set $r
= 3\Nobs$, noting that each sensor location provides three displacement measurements. In
practice, this choice of $r$ is conservative, and empiraclly we find that $r \approx
150$ is sufficient for the present problem.

Now, turning to the sensitivity analysis, we need to first parameterize the uncertainty
in the different auxiliary parameters consistently.  Suppose $\vartheta_i$, $i = 1,
\ldots, n$, are the auxiliary model parameters. Then, $\vartheta_i$ represents KLE
coefficients $\overline{\lambda}, \overline{\mu}, \lambda_i^{\rm KLE}, \mu_i^{\rm KLE}$
in the L\'ame parameters and coefficient $\nu_k, \nu_s $ of the Robin conditions.
Since we compare parameters with different physical units, we use the parameterization 
\begin{equation}
  \vartheta_i = (1 + \alpha \theta_i) \overline{\vartheta}_i,
  \quad i \in \{1, \ldots, n\},
\end{equation}
where $\overline{\vartheta}_i$'s are the nominal values, $\alpha > 0$ is a coefficient
modeling our uncertainty level in the nominal value, and $\theta_i$ is the perturbation
parameter we study in place of $\vartheta_i$.  Note that this parameterization assumes
the nominal parameter values are nonzero.  We assume $\theta_i \in [-1, 1]$ and let
$\alpha = 0.05$, i.e., we consider an uncertainty level of $5\%$ for each auxiliary
parameter.  The present parameterization also sets the stage for a global sensitivity
analysis, where $\theta_i$ will be assumed uniformly distributed random variables on the
interval $[-1, 1]$.

\begin{figure}[!htb]
  \centering
  \begin{subfigure}[b]{0.495\textwidth}
    \caption{
      Information Gain Sensitivities \vspace{.2em}
    }
    \includegraphics[width=\textwidth]{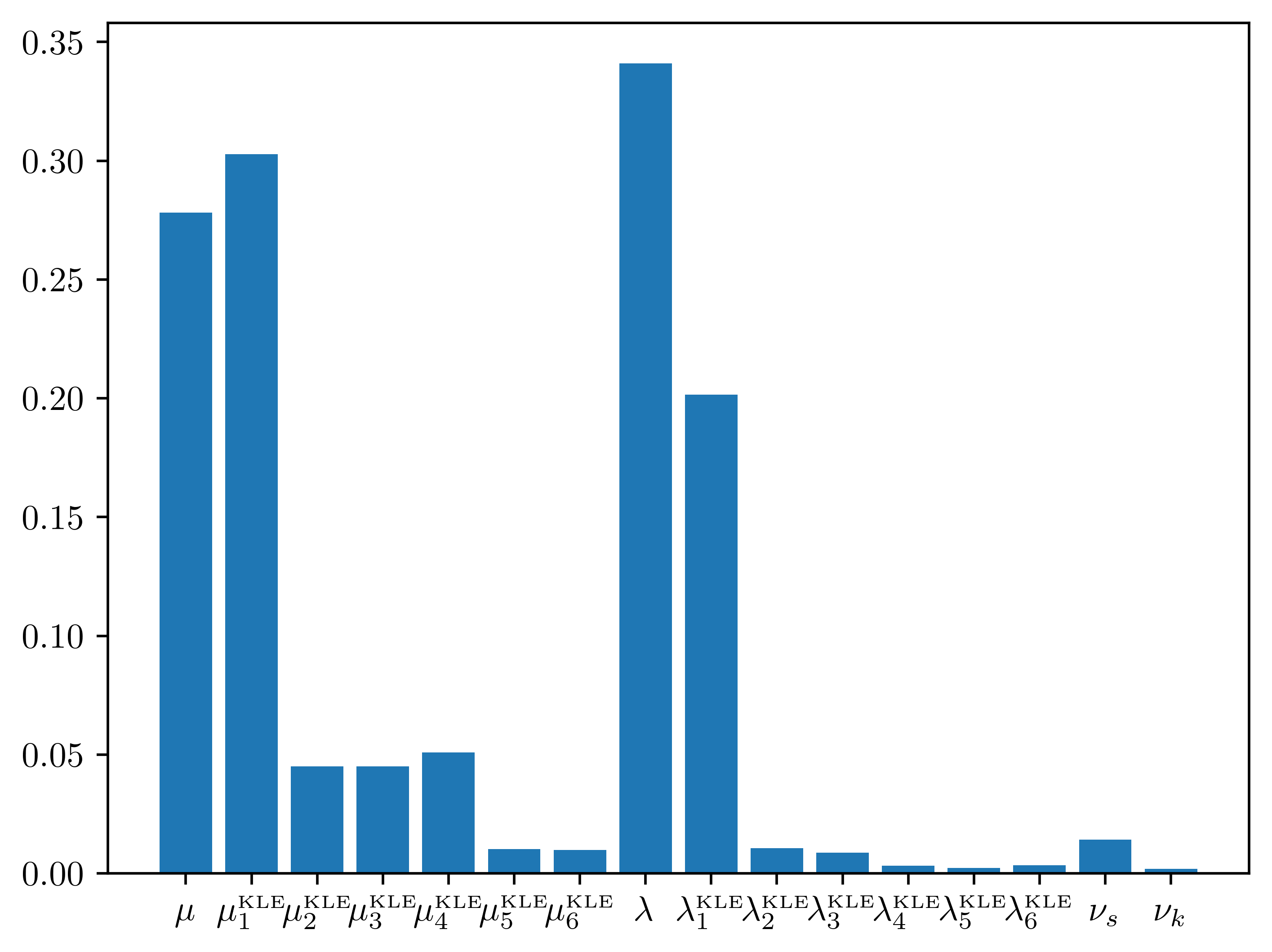}
  \end{subfigure}
  \begin{subfigure}[b]{0.495\textwidth}
    \caption{
      Expected Information Gain Sensitivities \vspace{.2em}
    }
    \includegraphics[width=\textwidth]{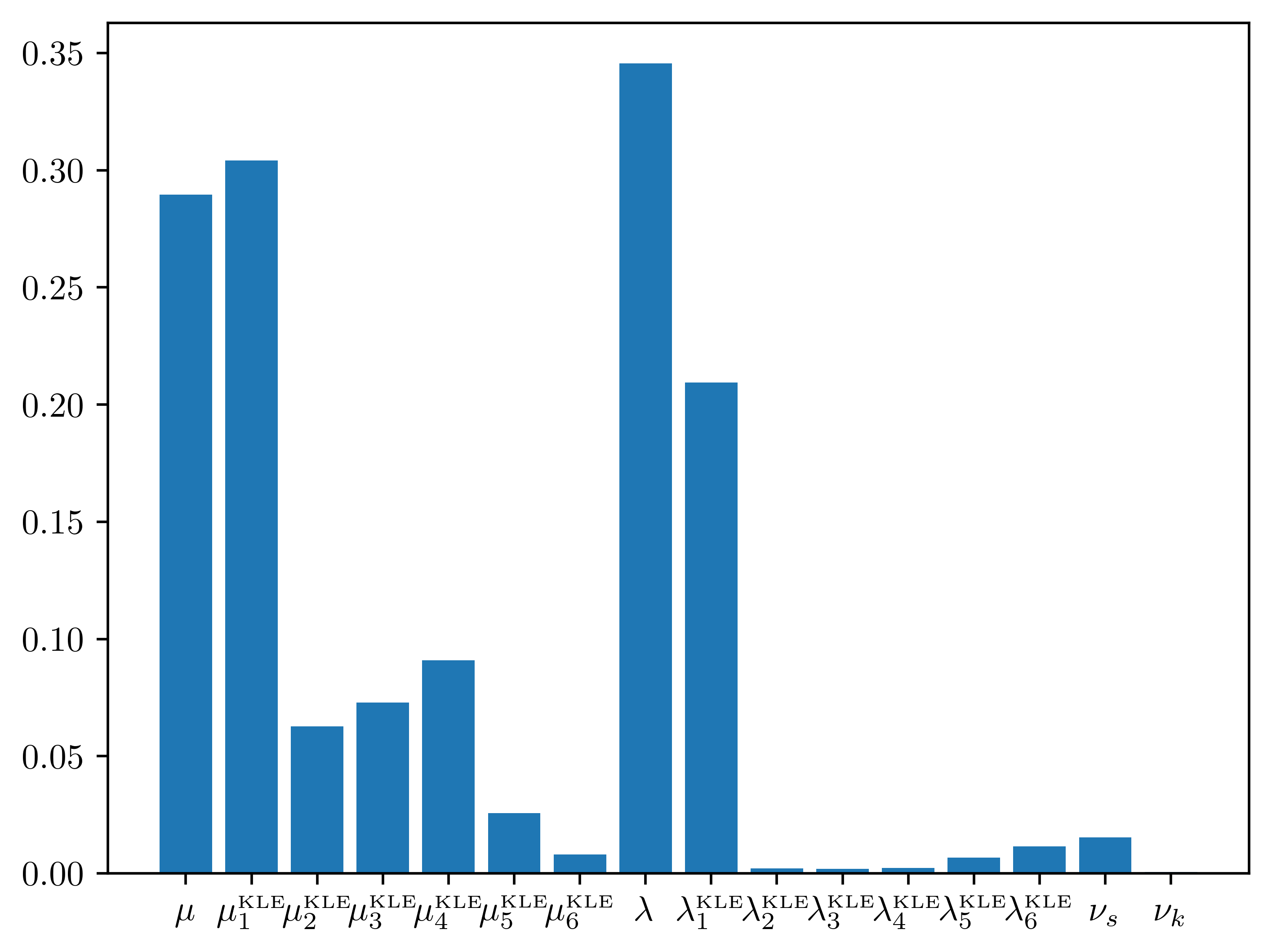}
  \end{subfigure}
  \caption{
    Sensitivity of the information gain and its expectation over data with respect to
    auxiliary parameters. Note, sensitivities are considered as absolute values of the 
    derivatives.
  } \label{fig:sensitivity}
\end{figure}
We use \Cref{alg:local_sensitivity} to compute the sensitivities of the information gain and its
expectation over data with respect to the auxiliary parameters. The results are reported
in \Cref{fig:sensitivity}.  We find that the L\'{a}me parameters hold significantly more
sensitivity than $\nu_k$ and $\nu_s$ from the boundary conditions. Moreover, for both
L\'{a}me parameters, the mean value and the first KLE mode are the most influential.
From an uncertainty quantification perspective, these results can be viewed as a
positive for the well-foundedness of the model. Indeed, both $\nu_k$ and $\nu_s$ are
artificial numerical parameters that are used due to the truncation of the domain.
Furthermore, the sensitivity being dominated by the first KLE modes indicates that, for
this model, it may be acceptable to ignore more oscillatory modes of the L\'{a}me
parameters as long as the mean is correctly chosen in the first place.

\begin{figure}[!htb]
  \centering
  \includegraphics[width=\textwidth]{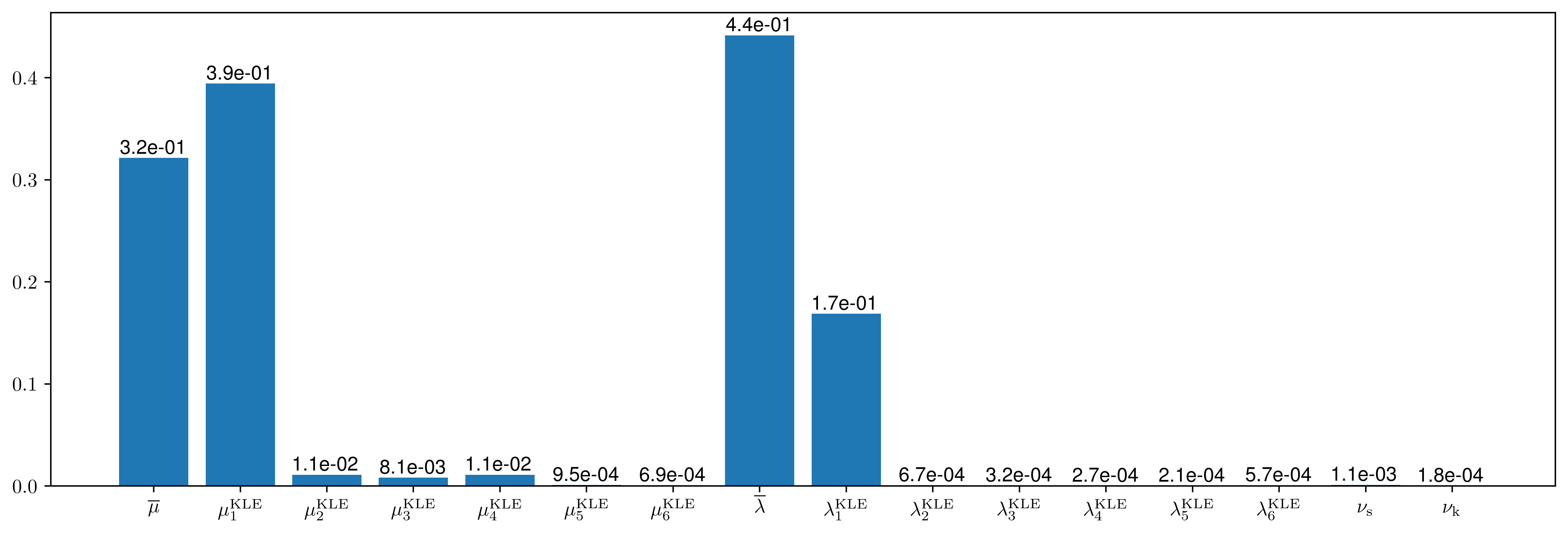}
  \caption{
    Derivative-based Sobol indices of the information gain with respect to the auxiliary
    parameters for the elasticity model. 
  } \label{fig:sobol}
\end{figure}
Finally, we study the global sensitivity of the information gain with respect to
the auxiliary parameters. As discussed in \Cref{sec:math-prelim-sobol}, this can
be done by considering the derivative-based upper bounds on the total Sobol
indices.  
The integrals in~\cref{eq:sobol-bound} are approximated via sample averaging, as
discussed in \cref{eq:variance-sample-average} and \cref{eq:dgsm-sample-average} and a
total of 500 samples were shared between the two estimators. We report the computed
Sobol index bounds in~\Cref{fig:sobol}. These results are generally consistent with the
behavior seen in the local sensitivity results in \Cref{fig:sensitivity}.  In
particular, both $\nu_k$ and $\nu_s$ coefficients and the higher order KLE modes, except
$\mu_4^\mathrm{KLE}$, are unimportant.

\section{Conclusion} \label{sec:conclusion} 
In this article, we have outlined a scalable approach for sensitivity analysis of the
information gain and expected information gain with respect to auxiliary model
parameters.  Our approach builds on low-rank spectral decompositions,  adjoint-based
eigenvalue sensitivity analysis, and concepts from post-optimal sensitivity analysis.
As seen in our computational examples, the (expected) information gain can exhibit
complex dependence on auxiliary model parameters.  This makes a practical sensitivity
analysis framework important. The present work makes a foundational step in this
direction. 

A natural next step is extending the presented sensitivity analysis framework to the
case of nonlinear Bayesian inverse problems. A possible approach is to utilize a Laplace
approximation to the posterior. The present approach can be adapted to that setting.
Though, it will be important to assess the extent to which the sensitivity of the
Laplace approximation to the posterior with respect to the auxiliary parameters is
indicative of the sensitivity of the true posterior distribution to these parameters.

The proposed sensitivity analysis framework is also synergistic to efforts in the area
of Bayesian optimal experimental design
(OED)~\cite{ChalonerVerdinelli95,Alexanderian_21} and OED under
uncertainty~\cite{KovalAlexanderianStadler20,AlexanderianPetraStadlerEtAl21,
FengMarzouk19,AlexanderianNicholsonPetra22}.  Specifically, for OED under uncertainty,
sensitivity analysis can reveal the most influential uncertain auxiliary parameters.
Then, the OED problem may be formulated by focusing only on the most important auxiliary
parameters.  This process can potentially reduce the computational burden of OED under
uncertainty substantially. Investigating such parameter dimension reduction approaches
in the context of OED under uncertainty is an interesting area of future work.

\section{Acknowledgment} \label{sec:ack}
The work of A.~Alexanderian and A.~Chowdhary was supported in part by US National
Science Foundation grant DMS-2111044. The authors would also like to thank Tucker
Hartland and Noemi Petra for their help in the slip inference computational example.

%\bibliography{References}
\printbibliography

\end{document}